\documentclass[titlepage,a4paper, 11pt]{article}
\usepackage[T1]{fontenc}
\usepackage[latin1]{inputenc}
\usepackage{amsmath}
\usepackage{amssymb}

\usepackage{amsthm}

\newtheorem{ex}{Example}[section]
\newtheorem{Def}{Definition}[section]
\newtheorem{Not}{Notational convention}[section]
\newtheorem{rem}{Remark}[section]
\newtheorem{Th}{Theorem}[section]
\newtheorem{lem}{Lemma}[section]

\newtheorem{cor}{Corollary}[section]

\newcommand{\PP}{\mathbb{P}}
\newcommand{\EE}{\mathbb{E}}
\newcommand{\II}{\mathbb{I}}
\newcommand{\CC}{\mathbb{C}}
\newcommand{\RR}{\mathbb{R}}
\newcommand{\QQ}{\mathbb{Q}}
\newcommand{\NN}{\mathbb{N}}

\newcommand{\cB}{{\cal B}}

\newcommand{\cC}{{\cal C}}
\newcommand{\cE}{{\cal E}}
\newcommand{\cI}{{\cal I}}

\newcommand{\cF}{{\cal F}}

\newcommand{\cG}{{\cal G}}

\newcommand{\cT}{{\cal T}}
\newcommand{{{\cadlag}}}{c\`adl\`ag}

\newcommand{\esssup}{\mathrm{ess \ sup} \ }

\begin{document}
\sloppy

\sffamily
\pagestyle{plain}

\renewcommand{\baselinestretch}{1}

\title{Continuity corrections for certain perpetual American and Bermudan options on multiple assets}
\author{Frederik Herzberg\thanks{Abteilung f\"ur Stochastik, Institut f\"ur Angewandte Mathematik, Universit\"at Bonn, D-53115 Bonn, Germany ({\tt herzberg@wiener.iam.uni-bonn.de})} \thanks{Mathematical Institute, University of Oxford, Oxford OX1 3LB, England}}

\date{}

\maketitle

\begin{abstract} 
In a general Feller martingale market with several assets, the existence of optimal exercise regions for multi-dimensional Bermudan options can be established by reference to Neveu's theory of Snell envelopes -- and also, as will be shown, more directly from standard martingale arguments combined with the strong Markov property. Based on this, in the framework of a log-L\'evy martingale market, explicit formulae and asymptotic results on the perpetual American-Bermudan (barrier-like) put option price difference (``continuity correction'') near the exercise boundary will be proven, under the -- of course, fairly restrictive -- assumption that the logarithmic optimal exercise region, subject to the barrier, does not depend on the time mesh size and is, up to translation, a half-space. 

For this sake, Wiener-Hopf type results by Feller will be generalized to higher dimensions. It will be shown that an extrapolation from the exact Bermudan prices to the American price cannot be polynomial in the exercise mesh size in the setting of many common market models, and more specific bounds on the natural scaling exponent of the non-polynomial extrapolation for a number of (both one- and multi-dimensional) market models will be deduced.
%
\end{abstract}

\tableofcontents
\pagebreak 

\noindent


\renewcommand{\baselinestretch}{1}

\subsection{Introduction}

{\em American options} are financial securities that can be exercised at any future time before maturity, or, in the case of perpetual options, at any date in the future. {\em Bermudan options}, on the other hand, can only be exercised at certain dates in the future (and therefore are sometimes also referred to as {\em discretely sampled American options}). An {\em option on a multiple asset} (or: {\em multi-dimensional option}) is a contract whose payoff function depends on more than one asset. An (American or Bermudan) barrier option can only be exercised if some asset price falls below or rises above a certain level. 

What the holder of an American or Bermudan option can gain from exercising the option depends on the exercise strategy she chooses. The issuer of the option has to set as price the least upper bound of all expected gains from such exercise strategies. 

In a memory-less market, this supremum coincides with the supremum over all exercise strategies that only look at the current price, and this supremum, of course, is a maximum, since at each exercise time one will only have to check whether the payoff one would get from exercising now is still less than the current option price (which in turn is taken to be the supremum mentioned above). By choosing smaller and smaller steps between the (discrete) exercise times of a Bermudan option, one can approximate an American price.

Whilst the reasoning of the previous paragraph has long been brought to mathematical precision for one-dimensional options, and nowadays seems to be fairly well-understood for higher dimensions as well, it appears still to be quite difficult to find references that rigorously prove the existence of exercise regions for multi-dimensional Americans or Bermudans. Therefore, this paper fills a gap in the existing literature by providing a rigorous introduction to (multi-dimensional) American and Bermudan option pricing in Section \ref{formalintro}.

We will thus first of all give a rigorous account of the existence of optimal exercise regions for Bermudan options on multiple assets, their coinciding with immediate exercise regions and (as an immediate Corollary) their time-stationarity in the case of perpetual options. Explicit characterisations of immediate exercise regions for (non-pereptual) Bermudan options can be found in the work of Broadie and Detemple \cite{BD}, Ekstr\"om \cite{Eks} as well as Peskir \cite{Pesk}, whereras the immediate exercise boundary of one-dimensional perpetual Bermudan options is discussed in an article by Boyarchenko and Levendorskii \cite{BL}. (In the spirit of Boyarchenko and Levendorskii's paper, we shall look at the operator equations corresponding to -- multi-dimensional -- perpetual Bermudan pricing problems in an Appendix.)

Having proven the existence and time-stationarity of the exercise regions for Bermudan options, we shall then move on to studying continuity corrections -- that is, the differences of an American price and a Bermudan price of a given exercise mesh size -- of a perpetual put option.

In particular, we shall derive asymptotic bounds on these corrections (conceived of as a function of the Bermudan's exercise mesh size) as the Bermudan exercise mesh size tends to zero. For this purpose, a Wiener-Hopf type result of Feller's \cite[p. 606, Lemma 3]{F} will be generalised to higher dimensions. 

The motivation for finding such continuity corrections is, of course, to be able to extrapolate from a finite number of Bermudan option prices to an approximation to the American price (both with the same particular payoff function). Therefore, the problem has attracted interest for quite some time, and computationally useful results for one-dimensional options have been established as well. One can find these scaling results in the works by Broadie, Glasserman and Kou \cite{BGK}, by Howison \cite{How} as well as by Howison and Steinberg \cite{HS}. 

The class of payoff functions $g$ to which our continuity correction results will be applicable is the class of half-space barrier-like payoff functions -- given a discount rate $r>0$, the {\em barrier-like payoff functions} are defined as exactly those payoff functions with the property that for sufficiently small excercise mesh sizes $s\leq s_0$, the optimal logarithmic exercise region $G:=G_s$ does not depend on $s$. If this $G$ is, up to translation, a set that is closed under addition in $\RR^d$ with its complement also being closed under addition, then the payoff function will be called a {\em half-space barrier-like payoff function}. We will also show that all convex, closed and $+$-closed sets with $0$ on the boundary whose complements are also $+$-closed, are merely half spaces in the following sense: If $H\subseteq \RR^d$ is closed, convex and $+$-closed with $\complement H$ also being $+$-closed, and $0\in\partial H$, then there is a $y_H\in\RR^d$ such that $H=\left\{x\in\RR^d\ : \ {^t}x\cdot y_H\geq 0\right\}$. 

Note that the price of a (possibly multidimensional) knock-in barrier option with payoff function $g$ where the possible exercise region $G$ determined by the barrier(s) is a subset of the optimal exercise region of the corresponding American option will be the same as the price of the American/Bermudan option with barrier-like payoff function $g\chi_G$, since the option will be exercised immediately after the barrier has been hit. This vindicates the term ``barrier-like payoff functions''.

This paper is based on Chapter 1 and Appendix A of the author's thesis \cite{H06}.

\subsection{Notation}

We are following largely standard probabilistic notation, as can be found for instance in the works by It\^o and McKean jr \cite{IK} or Revuz and Yor \cite{RY}. In particular, if $\PP$ is a probability measure on a $\sigma$-algebra $\cC$, $Z$ a random variable and $C\in\cC$, then $\EE\left[Z,C\right]=\int_C Z d\PP$, and if $\PP[C]>0$, then $\EE[Z|C]=\frac{1}{P}\EE\left[Z,C\right]$, whilst the conditional expectation of $Z$ with respect to a sub-$\sigma$-algebra $\cC_0\subseteq \cC$ will be denoted by $\EE\left[Z|\cC_0\right]$

Both $A\subseteq B$ and $A\subset B$ for sets $A$ and $B$ will mean that $A$ is a subset of $B$ (possibly $A=B$), whereas we shall write $A\subsetneq B$ to express that $A$ is a proper subset of $B$.

Finally, for any subset $A\subseteq \RR^d$, $\complement A$ shall denote its complement $\complement A:=\RR^d\setminus A\subseteq\RR^d$.

\section{Definitions and basic facts on Bermudan and American options}

\label{formalintro}

\subsection{Terminology}

Our first definition is a notational convention.

\begin{Def}\label{lieexp} Let $d\in\NN$. By $\exp:\RR^d\rightarrow {\RR_{>0}}^d$ and $\ln:{\RR_{>0}}^d\rightarrow\RR^d$ we denote componentwise exponentiation and taking natural logarithms componentwise, respectively.
\end{Def}
\begin{rem} For any $d\in\NN$, $\RR^d$ is a Lie group with respect to componentwise multiplication $\cdot:(x,y)\mapsto (x_iy_i)_{i\in\{1,\dots,d\}}$. Its Lie algebra is the vector space $\RR^d$ with its usual (componentwise) addition. The exponential map from the Lie algebra $(\RR^d,+)$ into the Lie group $(\RR^d,\cdot)$ is componentwise exponentiation $\exp:x\mapsto\left(e^{x_i}\right)_{i\in\{1,\dots,d\}}$. Therefore the abbreviation introduced in Definition \ref{lieexp} is consistent with standard notation.
\end{rem}

\begin{Def} Let $T$ be a positive real number. Consider a real-valued stochastic process $X:=(X_t)_{t\in[0,T]}$, adapted to a filtered probability space $\left(\Omega,(\cF_t)_{t\in[0,T]},P\right)$. We will call $X$ a {\em logarithmic price process for an asset with continuous dividend yield $\delta$} (for short, a {\em logarithmic price process} or simply {\em log-price process}), if and only if there exists a probability measure $Q$ equivalent to $P$ on $\cF_T$ and a constant $r>0$ such that the stochastic process $\exp\left(X_t-rt+\delta t\right)_{t\in[0,T]}$ is a martingale with respect to the filtration $\cF:=(\cF_t)_{t\in[0,T]}$ and the probability measure $Q$. In this case, such a $Q$ is called a {\em martingale measure} and $r$ a {\em market price of risk} or a {\em dicsount rate} for the stochastic process $X$ and the probability measure $P$.
\end{Def}

\begin{Def} Let $d\in\NN$. A {\em $d$-dimensional basket} is a $d$-tuple of logarithmic price processes such that there exists a probability measure $Q$ and a market price of risk $r>0$ such that $Q$ is a martingale measure and $r$ a market price of risk for all components of the $d$-tuple.
\end{Def}

For the rest of this section, we will adopt the terminology and the notation for Markov processes of Revuz and Yor \cite{RY}.

In particular, for all probability measures $\nu$ on $\cB\left(\RR^d\right)$, $Q_\nu$ is the probability measure induced by the transition function $\left(Q_s\right)_{s\geq 0}$ via the Ionescu-Tulcea-Kolmogorov projective limit construction, cf. Revuz and Yor \cite[Theorem 1.5]{RY}).

For any $d\in\NN$, we will denote the $\sigma$-algebra of Borel subsets of $\RR^d$ by $\cB\left(\RR^d\right)$.

\begin{Def} \label{Markovbasket}Let again $d\in\NN$. A family $Y:=\left(Y^x_\cdot\right)_{x\in\RR^d}$ of $\RR^d$-valued homogeneous Markov processes $Y^x$ adapted to a filtered probability space $\left(\Omega,\cF,\tilde Q\right)$ with respect to $\cF$, with transition function $\left(P_s\right)_{s\geq 0}$ and initial measure $\delta_x$, is called a {\em $d$-dimensional Markov basket with dividend yields $\delta_1,\dots,\delta_d>0$} if and only if there is a homogeneous transition function $\left(Q_s\right)_{s\geq 0}$ on the measurable space $\left(\RR^d,\cB\left(\RR^d\right)\right)$ and a constant $r>0$ such that the following three assertions hold:
\begin{enumerate}
\item The process $Y^x_\cdot$ is a Markov process with transition function $\left(Q_s\right)_{s\geq 0}$ with respect to $\cF$ for all $x\in\RR^d$.
\item The process $\left(\left(\exp\left(\left({Y^x}_t\right)_i-(r-\delta_i)t\right)\right)_{i\in\{1,\dots,d\}}\right)_{t\in[0,T]}$ is a martingale with respect to $\cF$ and $Q_{\delta_x}$.
\item The measures $P_{\delta_x}$ and $\PP^x:=Q_{\delta_x}$ are equivalent for all $x\in\RR^d$.
\end{enumerate}

In this case, $\left(\PP^x\right)_{x\in\RR^d}$ is called a {\em family of martingale ({\rm or:} risk-neutral) measures associated with $Y$}, and $r$ is called the {\em discount rate for $Y$}. 

The expectation operator for the probability measure $\PP^x$ will be denoted by $\EE^x$ for all $x\in\RR^d$.

If the transition function $P$ is a Feller semigroup, then we shall refer to $Y$ as a {\em Feller basket}.

If $P$ is a translation-invariant Feller semigroup, we shall call $Y$ a {\em L\'evy basket}.

\end{Def}

\begin{rem} A priori, it is not clear if there are logical connections between the three assertions in the previous Definition \ref{Markovbasket}, in particular the author does not know whether the third assertion implies the first one.
\end{rem}

\begin{Not} If no ambiguity can arise, we will drop the superscript of a Markov basket. Thus, in the notation of Definition \ref{Markovbasket}, we set $$\EE^x\left[\left.f\left(Y_{\tau_1},\dots,Y_{\tau_n}\right)\right|\cF_s\right]:=\EE^x\left[\left.f\left(Y_{\tau_1}^x,\dots,Y_{\tau_n}^x\right)\right|\cF_s\right]$$ for all $s\geq 0$, $n\in\NN$ and $n$-tuples of stopping times $\vec \tau=(\tau_1,\dots,\tau_n)$ whenever $f:\RR^n\rightarrow\RR$ is nonnegative or $f\left(Y_{\tau_1}^x,\dots,Y_{\tau_n}^x\right)\in L^1\left(\PP^x\right)$. Here we are using the term ``stopping time'' as a synonym for $\RR_+$-valued stopping time, that is a stopping time with values in $[0,+\infty]$. 

Also, since we are explicitly allowing stopping times (with respect to the filtration generated by a process $X$) to attain the value $+\infty$, we stipulate that the random variable $f\left(X_\tau\right)$ (for any Lebesgue-Borel measurable function $f$) should be understood to be multiplied by the characteristic function of the event $\left\{\tau<+\infty\right\}$. Formally, this can be done by introducing a constant $\Delta\not\in \RR^d$, called {\em cemetery}, and stipulating that $X_\tau=\Delta$ on $\left\{\tau=+\infty\right\}$ and $f(\Delta)=0$ for all measurable functions $f$ (cf. e.g. Revuz and Yor \cite[pp 84,102]{RY}).
\end{Not}

We will not formally define what we mean by an option itself, but rather provide definitions for the concepts of {\em expected payoffs} and {\em prices} for certain classes of financial derivatives.

\begin{Def}

Consider a $d$-dimensional Markov basket $Y$ with an associated family $\PP^\cdot$ of martingale measures and discount rate $r>0$. 

The expected payoff of a {\em Bermudan option with (log-price) payoff function $g:\RR^d\rightarrow \RR_{\geq 0}$ on the underlying Markov basket $Y$ with exercise times in $J\subset[0,+\infty)$, log start-price $x$ and maturity $T\in[0,+\infty]$} is defined to be $$U^J(T)(x):=U^J_g(T)(x):=\sup_{\tau\text{ stopping time}, \ \tau(\Omega)\subseteq J\cup\{+\infty\} }\EE^x\left[e^{-r(\tau\wedge T)}g\left(X_{\tau\wedge T}\right)\right].$$ 

The expected payoff of a {\em perpetual Bermudan option} is the expected payoff of a Bermudan option of maturity $+\infty$.

The expected payoff of a {\em Bermudan option with exercise mesh size $h>0$} is the expected payoff of a Bermudan option with exercise times in $h\cdot\NN_0$ .

The expected payoff of an {\em American option} is the expected payoff of a Bermudan option with exercise times in $[0,+\infty)$.

We shall call the expected payoff of a Bermudan option (or an American option) a {\em Bermudan option price} (or an {\em American option price}) if and only if the martingale measures associated with the underlying basket are unique (that is, if the market model described by $P$, $\cF$ and $X$ is {\em complete}).

\end{Def}

In recent years, there has been increasing interest in incomplete market models that are governed by general L\'evy processes as log-price processes, as is not only witnessed by a tendency in research papers to focus on L\'evy process settings (for instance Boyarchenko and Levendorskii \cite{BL}; Asmussen, Avram and Pistorius \cite{AAP}; \O ksendal and Proske \cite{OP}, to take a random sample). Even textbooks, such as Karatzas' \cite{K} and Mel'nikov's \cite{M} introductory works, are putting considerable emphasis on incomplete markets. Finally, ``L\'evy finance'' has already been treated in survey articles intended for a general mathematical audience, e.g. Applebaum's article \cite{A04}. 

Whilst quite a few of our results will apply only to the Black-Scholes model, some of our arguments also work for market models where the logarithmic price process merely needs to be a L\'evy process. However, the theory of incomplete markets (and all L\'evy models other than the Black-Scholes model are incomplete) is not the focus of this thesis, therefore we did not endeavour to go beyond the Black-Scholes model where this caused technical difficulties rather than making proofs easier to read. 

Whilst there are some points to be made about market failures on stock markets that might entail arbitrage opportunities (for example, when assets are traded simultaneously on several stock exchanges, or in the event of insider trading), the transaction costs to exploit these arbitrage opportunities usually tend to be close to the actual gain that can be achieved through taking advantage of the arbitrage. Therefore we shall, for the sake of mathematical simplicity, merely refer to the works of Corcos et al. \cite{CEMMS}, Imkeller \cite{I} as well as Imkeller, Pontier and Weisz \cite{IPW}, and impose a strict no-arbitrage assumption -- which under certain regularity conditions on the basket is equivalent to the existence of an equivalent martingale measure (a measure equivalent to the market model under which the discounted asset prices are martingales), cf. Karatzas \cite[Theorem 0.2.4]{K} and references therein.

\begin{ex}[A few common examples] 
\begin{enumerate} 
\item The price of a European call option on a single asset with maturity $T$ and strike price $K$ is the price of a Bermudan option with the set of exercise times being the singleton $\{T\}$ and the (log-price) payoff function $\left(\exp(\cdot)-K\right)\vee 0$.
\item The price of a perpetual American put of exercise mesh size $h>0$ on the arithmetic average of two assets in an underlying basket with strike price $K$ is the price of a Bermudan put option with the set of exercise times being the whole of the half-line $[0,+\infty)$, the maturity being $T$ and the payoff function $\left(K-\frac{\exp\left((\cdot)_1\right)+\exp\left((\cdot)_2\right)}{2}\right)\vee 0$.
\item Consider a perpetual Bermudan call option on a single asset that continuously pays dividends at a rate $\delta$ and whose logarithm follows a Markov process $Z$ adapted to some probability space $\left(\Omega,(\cF_t)_{t\geq 0},P\right)$. Then, in order to exclude arbitrage, we will have to require the existence of a family of measures $\PP^\cdot$ such that each $\PP^x$ is equivalent to $P^x$ (in particular, $\PP^{x}_{Z_0}=\delta_x$) and such that $\left(e^{-rt+\delta t+ Z_t}\right)_{t\geq 0}$ is a $\PP^x$-martingale for all $x\in \RR^d$. The expected payoff of the option will then be $$\tilde U^{h\cdot \NN_0}(\cdot)=\sup_{\tau\text{ stopping time}, \ \tau(\Omega)\subseteq h\NN_0\cup\{+\infty\} }\EE^\cdot\left[e^{-r\tau}\left(e^{Z_{\tau}}-K\right)\vee 0\right]$$
\end{enumerate}
\end{ex}

\subsection{Convergence of Bermudan to American prices}

In this paragraph, we shall give a formal proof in a general setting that the price of a Bermudan option with equidistant exercise times (of mesh $h$) before maturity indeed converges to the American price as $h$ tends to $0$.

As an auxiliary Lemma, let us remark the following elementary observation:

\begin{lem}[Lower semi-continuity of $\sup$] If $I$ is a set and $\left(a_{k,\ell}\right)_{\ell\in I,k\in\NN_0}$ is a family of real numbers, then $$\sup_{\ell\in I} \liminf_{k\rightarrow\infty} a_{k,\ell}\leq \liminf_{k\rightarrow\infty} \sup_{\ell\in I} a_{k,\ell}.$$
\end{lem}
\begin{proof} We have the trivial estimate $$\sup_{\ell}a_{k,\ell}\geq a_{k,\ell_0}$$ for all $k\in\NN_0$ and $\ell_0\in I$, therefore for all $n\in\NN$ and $\ell_0\in I$, $$\inf_{k\geq n }\sup_{\ell}a_{k,\ell}\geq \inf_{k\geq n} a_{k,\ell_0},$$ thus $$\inf_{k\geq n }\sup_{\ell}a_{k,\ell}\geq \sup_{\ell_0}\inf_{k\geq n} a_{k,\ell_0},$$ hence $$\sup_n\inf_{k\geq n }\sup_{\ell}a_{k,\ell}\geq \sup_n\sup_{\ell_0}\inf_{k\geq n} a_{k,\ell_0}=\sup_{\ell_0}\sup_n\inf_{k\geq n} a_{k,\ell_0}.$$ This is the assertion.
\end{proof}

This estimate enables us to prove the following Lemma that is asserting the approximability of expected payoffs or prices of American options by sequences of expected payoffs or prices of Bermudan options, respectively.

\begin{lem}\label{dyadic limits generic}
Let $d\in\NN$, $T>0$ ($T=+\infty$ also possible), $x\in\RR^d$, and let $X$ be a $d$-dimensional Markov basket. Suppose the payoff function $g\geq 0$ is bounded and continuous, and assume $X$ has a modification with continuous paths. 

Then, if the American expected payoff $U^{[0,+\infty)}(T)(x)$ is finite, one has \begin{eqnarray*}U^{[0,+\infty)}(T)(x)& =& \lim_{h\downarrow 0}U^{h\NN_0}(T)(x) \\ &=& \sup_{k\in\NN} U^{2^{-k}\NN_0}(T)(x).\end{eqnarray*}
\end{lem}

\begin{proof} Consider a sequence $(h_k)_{k\in\NN_0}\in\left(\RR_{>0}\right)^{\NN_0}$ such that $h_k\downarrow 0$ as $k\rightarrow\infty$. Choose a sequence of stopping times $(\tau_\ell)_{\ell\in\NN_0}$ such that for all $x\in\RR^d$, $$\sup_{\tau\text{ stopping time}}\EE^x\left[e^{-r(\tau\wedge T)}g\left(X_{\tau\wedge T}\right)\right]=\sup_{\ell}\EE^x\left[e^{-r(\tau_\ell\wedge T)}g\left(X_{\tau_\ell\wedge T}\right)\right]$$ and define $$\tau_{\ell,k}:= \inf\left\{t\in h_k\NN_0 \ : \ t\geq \tau_\ell\right\}.$$ Then, due to the continuity conditions we have imposed on $g$ and on the paths of (a modification of) the basket $X$, we get $$\sup_{\ell} e^{-r\left(\tau_\ell\wedge T\right)}g\left(X_{\tau_\ell\wedge T}\right)= \sup_{\ell} \lim_{k\rightarrow\infty} e^{-r\left(\tau_{\ell,k}\wedge T\right)}g\left(X_{\tau_{\ell,k}\wedge T}\right)$$ and hence by the lower semi-continuity of $\sup$, one obtains $$\sup_{\ell} e^{-r\left(\tau_\ell\wedge T\right)}g\left(X_{\tau_\ell\wedge T}\right)\leq \liminf_{k\rightarrow\infty} \sup_{\ell} e^{-r\left(\tau_{\ell,k}\wedge T\right)}g\left(X_{\tau_{\ell,k}\wedge T}\right).$$ Now we can use the Montone Convergence Theorem and Lebesgue's Dominated Convergence Theorem (this is applicable because of the boundedness of $g$) to swap limits/suprema with the expectation operator. Combining this with the specific choice of the sequence $\left(\tau_\ell\right)_{\ell\in\NN_0}$, this yields for all $x\in\RR^d$, \begin{eqnarray*} &&\sup_{\tau\text{ stopping time}}\EE^x\left[e^{-r(\tau\wedge T)}g\left(X_{\tau\wedge T}\right)\right]=\sup_{\ell}\EE^x\left[e^{-r(\tau_\ell\wedge T)}g\left(X_{\tau_\ell\wedge T}\right)\right] \\ &=& \EE^x\left[\sup_{\ell}e^{-r(\tau_\ell\wedge T)}g\left(X_{\tau_\ell\wedge T}\right)\right] \leq \EE^x\left[\liminf_{k\rightarrow\infty} \sup_{\ell} e^{-r\left(\tau_{\ell,k}\wedge T\right)}g\left(X_{\tau_{\ell,k}\wedge T}\right)\right]\\ &\leq& \liminf_{k\rightarrow\infty}\EE^x\left[ \sup_{\ell} e^{-r\left(\tau_{\ell,k}\wedge T\right)}g\left(X_{\tau_{\ell,k}\wedge T}\right)\right] \\ &=& \liminf_{k\rightarrow\infty}\sup_{\ell} \EE^x\left[ e^{-r\left(\tau_{\ell,k}\wedge T\right)}g\left(X_{\tau_{\ell,k}\wedge T}\right)\right] \\&\leq& \liminf_{k\rightarrow\infty}\sup_{\tau(\Omega)\subseteq h_k\NN_0\cup\{+\infty\}} \EE^x\left[ e^{-r\left(\tau\wedge T\right)}g\left(X_{\tau\wedge T}\right)\right] \\&\leq& \limsup_{k\rightarrow\infty}\sup_{\tau(\Omega)\subseteq h_k\NN_0\cup\{+\infty\}} \EE^x\left[ e^{-r\left(\tau\wedge T\right)}g\left(X_{\tau\wedge T}\right)\right] \\&\leq& \sup_{\tau\text{ stopping time}}\EE^x\left[e^{-r(\tau\wedge T)}g\left(X_{\tau\wedge T}\right)\right]. \end{eqnarray*} This finally gives \begin{eqnarray*}&&\sup_{\tau\text{ stopping time}}\EE^x\left[e^{-r(\tau\wedge T)}g\left(X_{\tau\wedge T}\right)\right]\\&=&\lim_{k\rightarrow\infty}\sup_{\tau(\Omega)\subseteq h_k\NN_0\cup\{+\infty\}} \EE^x\left[ e^{-r\left(\tau\wedge T\right)}g\left(X_{\tau\wedge T}\right)\right]=U^{h_k\NN_0}(T)(x). \end{eqnarray*} Since the left hand side does not depend on $(h_k)_k$, we conclude that $\lim_{h\downarrow 0}U^{h\NN_0}(T)(x)$ exists and is equal to $\sup_{\tau\text{ stopping time}}\EE^x\left[e^{-r(\tau\wedge T)}g\left(X_{\tau\wedge T}\right)\right]$.
\end{proof}

\section{Exercise regions}

An American/Bermudan option price coincides with the payoff that is expected if one exercises at the first possible entry of the log-price process into the immediate exercise region (which is a subset $G\subset \RR^d$ in case of a perpetual option and a subset $\cG\subset \RR^d\times [0,+\infty)$ for a non-perpetual). The immediate exercise region for a Markov basket with payoff function $g$ is defined as $$F^{J,T}:=\left\{(x,t)\in\RR^d\times J \ : \ U^J(T-t)(x)\leq g(x)\right\} $$ for a non-perpetual option with maturity $T$ and a set of exercise times $J$, and as $$F^h:=\left\{x\in\RR^d \ : \ U^{h\NN_0}(+\infty)(x)\leq g(x)\right\} $$ for perpetual options with exercise mesh size $h$; its optimality for a large class of options was established in the theory of optimal stopping (cf. e.g. Neveu \cite[Proposition VI-2-8]{Nev} or El Karoui \cite[Th\'eor\`eme 2.31]{ElK}, whose results need to be applied to the corresponding space-time Markov process) using the so-called {\em Snell envelope} (cf. Griffeath and Snell \cite{Snell}). In this Section, in addition to setting up notation, we shall give alternative elementary proofs for the optimality of immediate exercise regions for put options with a discrete set of exercise times. Explicit characterisations of immediate exercise regions for certain special cases have been proven in recent years: by Broadie and Detemple \cite{BD}, Paulsen \cite{P}, Ekstr\"om \cite{Eks} as well as Peskir \cite{Pesk}. 

For the rest of this Section, we would like to restrict our attention to hitting times where, no matter if the paths of the process are right-continuous and the target region closed, the infimum in the definition of a hitting time is always attained. We will do this by imposing the condition that the range of the stopping time be {\em discrete} in the following sense.

\begin{Def} A subset $I$ of a topological space $(X,\cT)$ is called {\em discrete} (with respect to $\cT$) if for all $x\in I$ there exists an open set $U\ni x$ such that $I\cap U=\{x\}$.
\end{Def}

Given any discrete subset $\cI$ of a $\RR_+$, it is, by density of the rationals in $\RR$, possible to find an embedding of $\cI$ into $\QQ_+$, thus $\cI$ must be countable. By an analogous argument, all discrete subsets of separable metric spaces must be countable.

A subset of a discrete set (with respect to a topology $\cT$) is again discrete with respect to the same topology $\cT$, and if $\cT$ comes from a linear order, the infimum of any discrete set is attained and therefore by definition a minimum.

\begin{Def} Given a discrete subset $\cI\subset [0,+\infty)$ and a Lebesgue-Borel measurable set $G\subset\RR^d$, often referred to as {\em exercise region}, we define the stopping time $$\tau_G^\cI:=\min\left\{t\in \cI \ : \ X_t\in G \right\},$$ (the superscript will be dropped when no ambiguity can arise) which is just the first (nonnegative) entry time in $\cI$ into $G$. If $\cG$ is a subset of space-time, that is $\cG\subset\RR^d\times[0,+\infty)$ rather than space (ie $\RR^d$) itself, we use the space-time process rather than just the process itself to give an analogous definition: $$\tau_\cG^\cI:=\min\left\{t\in t_0+\cI\ : \ (X_t,t)\in \cG \right\},$$ where $t_0$ is the time-coordinate at which the space-time process was started. Also, for $h>0$ we set $$\tau^h_G:=\tau_G^{h\NN},\quad \tau^h_\cG:=\tau_\cG^{h\NN}$$ to denote the first positive entry time in $h\NN_0$ into $G$ or $\cG$, respectively, whilst finally $\bar\tau_\cG^h:=\tau_\cG^{h\NN_0}$ and $\bar\tau_G^h:=\tau_G^{h\NN_0}$ will denotes the first nonnegative entry time into $\cG$ and $G$, respectively.

\end{Def}

For convenience, we will also adopt the following convention for this section:

\begin{Def} Let $\cI\subset[0,+\infty)$. A stopping time $\tau$ is called {\em $\cI$-valued} if the range of $\tau$, denoted by $\mathrm{ran \ } \tau$, is a subset of $\cI\cup\{+\infty\}$.
\end{Def}

The following Lemma, as well as its Corollary can be proven easily by resorting to the well-understood theory of optimal stopping and Snell envelopes, cf. e.g. Neveu and El Karoui. Our proof will be elementary.

\begin{lem}\label{exerciseboundaryLemma} Consider a discrete subset $\cI\subset [0,+\infty)$. Let $X$ be a $d$-dimensional basket with an associated risk-neutral measure $\PP$ and discount rate $r>0$. Suppose $g= (K-f )\vee 0$ and $e^{-r\cdot}f(X_{\cdot})$ is a $\PP$-submartingale. For all $\cI$-valued and $\PP$-almost surely finite stopping times $\tau$ there is a space-time region $B=\bigcup_{u\in \cI} \{u\}\times B_u$ such that $$\EE\left[e^{-r\left(\tau\wedge T\right)}g(X_{\tau\wedge T})\right]\leq \EE\left[e^{-r\left(\tau_B^\cI\wedge T\right)}g\left(X_{\tau_B^\cI\wedge T}\right)\right]$$
for all $T\in[0,+\infty)$ where $$\tilde\tau:=\tau_B^\cI=\inf\left\{u\in \cI\ : \ X_u\in B_u\right\}$$ and $\tau\geq \tau_B^\cI$ $\PP$-almost surely. If the set $\left\{ e^{-r \upsilon}g\left(X_{\upsilon}\right) \ : \ \upsilon\text{ $\cI$-valued stopping time} \right\}$ is uniformly $\PP$-integrable, then the latter inequality will also hold for $T=+\infty$. 

The Lemma holds in particular for $\cI=s\NN_0$ for arbitrary $s>0$.
\end{lem}
The condition of $\left\{ e^{-r \upsilon}g\left(X_{\upsilon}\right) \ : \ \upsilon\text{ $\cI$-valued stopping time} \right\}$ being uniformly $\PP$-integrable is what is known in Neveu's terminology \cite[e.g. Proposition 2.29]{Nev} as $\left( e^{-r \upsilon}g\left(X_{\upsilon}\right) \right)_{\upsilon\text{ $\cI$-valued stopping time} }$ being of class $(D)$.
\begin{proof}[Proof of Lemma \ref{exerciseboundaryLemma}] Firstly, we will treat the case of $T<+\infty$. Define $$\forall t\in\cI\quad B_t:=X_{\tau}\left(\left\{\tau=t\right\}\right)\subset \RR^d.$$ Let us first of all assume that \begin{equation}\label{gXtau>0onAt}\forall t\in \cI\quad\{\tilde\tau=t\}\cap\left\{\tau>T\right\}=\emptyset.\end{equation} and let us also for the moment suppose \begin{equation}\label{gXtau>0onCt}\forall t\in \cI\quad g(X_t)>0 \text{ a.s. on }\left\{X_t\in B_t\right\}.\end{equation} Both of these assumptions will be dropped at the end of the proof for the case $T<+\infty$ in order to show the Lemma in its full strength. Now, from equations (\ref{gXtau>0onAt}) and (\ref{gXtau>0onCt}) one may derive \begin{eqnarray} &&\EE\left[e^{-r\left(\tau\wedge T\right)}g(X_{\tau\wedge T}),\left\{\tilde \tau=t\right\}\right] \nonumber \\ &=& \EE\left[e^{-r\tau}g\left(X_{\tau}\right),\left\{\tilde \tau=t\right\}\right] \nonumber \\ &=& \EE\left[e^{-r\tau}\left(K-f(X_{\tau})\right),\left\{\tilde \tau=t\right\}\right] \nonumber\\ &=& \EE\left[e^{-r\left(\tau\wedge T\right)}\left(K-f(X_{\tau\wedge T})\right),\left\{\tilde \tau=t\right\}\right] \label{gbyK-f}\end{eqnarray} for all $t\in[0,T]\cap \cI$.

Furthermore, observe that $\tau\geq \tilde\tau $ a.s. Using Doob's Optional Stopping Theorem (see e.g. Varadhan \cite[Theorem 5.11]{Va}), we infer from our assumption of $e^{-r\cdot}f(X_{\cdot})$ being a $\PP$-submartingale with respect to the canonical filtration $\cF$ the assertion that $\left(e^{-r\upsilon}f\left(X_\upsilon\right)\right)_{\upsilon\in \left\{\tilde\tau\wedge T, \tau\wedge T\right\}}$ is a $\PP$-submartingale with respect to the filtration $\left\{\cF_{\tilde\tau\wedge T},\cF_{\tau\wedge T}\right\}$. Hence, if we combine this with equation (\ref{gbyK-f}) and note that $\left\{\tilde \tau=t\right\}=\left\{\tilde \tau\wedge T=t\right\}\in\cF_{\tau\wedge T}$ for all $t\in[0,T)\cap \cI$, we obtain for every $t\in[0,T)\cap \cI$, \begin{eqnarray}\nonumber &&\EE\left[e^{-r\left(\tau\wedge T\right)}g(X_{\tau\wedge T}),\left\{\tilde \tau=t\right\}\right]\\ \nonumber &= & \EE\left[e^{-r\left(\tau\wedge T\right)}\left(K-f(X_{\tau\wedge T})\right),\left\{\tilde \tau=t\right\}\right]\\ \nonumber &\leq & K\cdot \EE\left[e^{-r\left(\tau\wedge T\right)},\left\{\tilde \tau=t\right\}\right]-\EE\left[e^{-r\left(\tilde\tau\wedge T\right)}f\left(X_{\tilde\tau\wedge T})\right),\left\{\tilde \tau=t\right\}\right]\\\nonumber &\leq & K\cdot \EE\left[e^{-r\left(\tilde\tau\wedge T\right)},\left\{\tilde \tau=t\right\}\right]-\EE\left[e^{-r\left(\tilde\tau\wedge T\right)}f\left(X_{\tilde\tau\wedge T})\right),\left\{\tilde \tau=t\right\}\right] \\ &= & \EE\left[e^{-r\left(\tilde\tau\wedge T\right)}g\left(X_{\tilde\tau\wedge T})\right),\left\{\tilde \tau=t\right\}\right] .\end{eqnarray} On the other hand, since $\tilde \tau\leq \tau$, if $\tilde\tau\geq T$, then also $\tau\geq T$, entailing $$\tilde \tau \wedge T=T=\tau\wedge T \text{ on } \left\{\tilde\tau \geq T\right\}.$$ Summarising these last two remarks, one concludes \begin{eqnarray} \nonumber&&\EE\left[e^{-r\left(\tau\wedge T\right)}g\left(X_{\tau\wedge T}\right)\right]\\ \nonumber &=&\sum_{t\in \cI}\EE\left[e^{-r\left(\tau\wedge T\right)}g\left(X_{\tau\wedge T}\right),\left\{\tilde \tau=t\right\}\right] \\ \nonumber &=& \sum_{t\in \cI\cap [0,T)}\EE\left[e^{-r\left(\tau\wedge T\right)}g\left(X_{\tau\wedge T}\right),\left\{\tilde \tau=t\right\}\right] \\ \nonumber &&+\sum_{t\in \cI\cap [T,+\infty)}\EE\left[e^{-r\left(\tau\wedge T\right)}g\left(X_{\tau\wedge T}\right),\left\{\tilde \tau=t\right\}\right] \\ \nonumber &\leq & \sum_{t\in \cI\cap [0,T)}\EE\left[e^{-r\left(\tilde\tau\wedge T\right)}g\left(X_{\tilde\tau\wedge T}\right),\left\{\tilde \tau=t\right\}\right] \\ \nonumber && + \sum_{t\in \cI\cap [T,+\infty)}\EE\left[e^{-r\left(\tilde\tau\wedge T\right)}g\left(X_{\tilde\tau\wedge T}\right),\left\{\tilde \tau=t\right\}\right]\\ \nonumber &=& \sum_{t\in \cI}\EE\left[e^{-r\left(\tilde\tau\wedge T\right)}g\left(X_{\tilde\tau\wedge T}\right),\left\{\tilde \tau=t\right\}\right] \\ \label{boundary-submartingale} &=& \EE\left[e^{-r\left(\tilde\tau\wedge T\right)}g\left(X_{\tilde\tau\wedge T}\right)\right].\end{eqnarray} 

In order to complete the proof for the case of $T<+\infty$, let us show that the assumptions (\ref{gXtau>0onCt}) and (\ref{gXtau>0onAt}) are dispensable. 

We first of all simply define the stopping time $$\tau':= \chi_{\complement\left(\bigcup_{t\in \cI}\left\{\tau =t\right\}\cap\left\{g(X_t)>0\right\}\right)}\cdot\infty+\sum_{t\in \cI}\chi_{\left\{\tau=t\right\}\cap\left\{g(X_t)>0\right\}}\cdot t$$ and based on this definition, we would set $$\forall t\in\cI\quad B_t':=X_{\tau'}\left(\left\{\tau'=t\right\}\right).$$ Then $$B_t'=X_{t}\left(\left\{\tau'=t\right\}\right)\subset X_t\left(\left\{g(X_t)>0\right\}\right)\subset \left\{g>0\right\},$$ hence \begin{equation}\label{CforB'}\forall t\in \cI\quad g(X_t)>0 \text{ a.s. on }\left\{X_t\in B_t'\right\}.\end{equation} However, in any case $$\EE\left[ e^{-r\left(\tau'\wedge T\right)}g(X_{\tau'\wedge T})\right]=\EE\left[ e^{-r\left(\tau\wedge T\right)}g(X_{\tau\wedge T})\right].$$ 

Now, suppose at least one of the conditions (\ref{gXtau>0onAt}) and (\ref{gXtau>0onCt}) was not satisfied (if \eqref{gXtau>0onCt} holds, one may even replace $\tau'$ by $\tau$ in what follows). In this situation we consider the stopping time $$\tau''=\chi_{\left\{\tau'\leq T\right\}\cup\left(\left\{\tau'>T\right\}\cap \left\{\tau_{B'}^\cI\geq T\right\}\right)}\cdot\tau'+\chi_{\complement\left(\left\{\tau'\leq T\right\}\cup\left(\left\{\tau'>T\right\}\cap \left\{\tau_{B'}^\cI\geq T\right\}\right)\right)}\cdot\infty.$$ If one now defines $$\forall t\in\cI\quad B_t'':=X_{\tau''}\left(\left\{\tau''=t\right\}\right)$$ then 
$$\forall t\in\cI\cap [0,T]\quad B_t''=X_{\tau''}\left(\left\{\tau'=\tau''=t\right\}\right)=X_{\tau'}\left(\left\{\tau'=t\right\}\right)= B_t'$$ and \begin{eqnarray*}\forall t\in\cI\cap (T,+\infty)\quad B_t''&=&X_{\tau''}\left(\left\{\tau''=t\right\}\right)=X_{\tau''}\left(\left\{\tau''=\tau'=t\right\}\cap \left\{\tau_{B'}^\cI\geq T\right\}\right)\\ &\subset& X_{\tau'}\left(\left\{\tau'=t\right\}\right)= B_t',\end{eqnarray*} hence \begin{equation}\label{CforB''}\forall t\in \cI\quad g(X_t)>0 \text{ a.s. on }\left\{X_t\in B_t''\right\}\end{equation} (because of \eqref{CforB'} and we have just seen $B''_t\subset B_t'$ for all $t\in\cI$). Furthermore,
\begin{eqnarray*}\tau_{B''}^\cI &=&\chi_{\left\{\tau'\leq T\right\}\cup\left(\left\{\tau'>T\right\}\cap \left\{\tau_{B'}^\cI\geq T\right\}\right)}\cdot\tau_{B'}^\cI+\chi_{\complement\left(\left\{\tau'\leq T\right\}\cup\left(\left\{\tau'>T\right\}\cap \left\{\tau_{B'}^\cI\geq T\right\}\right)\right)}\cdot\infty\\ &=& \chi_{\left\{\tau''\leq T\right\}}\cdot \tau_{B'}^\cI + \chi_{\left\{+\infty>\tau''>T\right\}\cap \left\{\tau_{B'}^\cI\geq T\right\}}\cdot\tau_{B'}^\cI\\ && +\chi_{\complement\left(\left\{\tau'\leq T\right\}\cup\left(\left\{\tau'>T\right\}\cap \left\{\tau_{B'}^\cI\geq T\right\}\right)\right)}\cdot\infty\end{eqnarray*} (the first line because of $X_{\tau''}\left(\complement\left(\left\{\tau'\leq T\right\}\cup\left(\left\{\tau'>T\right\}\cap \left\{\tau_{B'}^\cI\geq T\right\}\right)\right)\right)=\emptyset$). Therefore $$\tau_{B''}^\cI\geq T\quad \text{ on }\left\{\tau''>T\right\},$$ thus (\ref{gXtau>0onAt}) holds for $\tau_{B''}^\cI$ instead of $\tilde\tau$ and $\tau''$ instead of $\tau$. But we have already proven (\ref{CforB''}). Therefore, analogously to the derivation of \eqref{boundary-submartingale} under the assumptions of both (\ref{gXtau>0onCt}) and (\ref{gXtau>0onAt}), we get $$\EE\left[ e^{-r\left(\tau''\wedge T\right)}g(X_{\tau''\wedge T})\right]\leq \EE \left[e^{-r\left(\tau_{B''}^\cI\wedge T\right)}g(X_{\tau_{B''}^\cI\wedge T})\right].$$ On the other hand, however, $$\EE\left[ e^{-r\left(\tau''\wedge T\right)}g(X_{\tau''\wedge T})\right]=\EE\left[ e^{-r\left(\tau'\wedge T\right)}g(X_{\tau'\wedge T})\right]$$ (as $\tau''=\tau'$ on $\{\tau< T\}$, as well as $\tau''\geq T$ on $\{\tau\geq T\}$, thus $\tau''\wedge T=T=\tau\wedge T $ on $\{\tau\geq T\}$) and we have already seen that $$\EE\left[ e^{-r\left(\tau'\wedge T\right)}g(X_{\tau'\wedge T})\right]=\EE\left[ e^{-r\left(\tau\wedge T\right)}g(X_{\tau\wedge T})\right].$$ Finally, $$\EE\left[ e^{-r\left(\tau\wedge T\right)}g(X_{\tau\wedge T})\right]\leq \EE \left[e^{-r\left(\tau_{B''}^\cI\wedge T\right)}g(X_{\tau_{B''}^\cI\wedge T})\right]$$ whence with $B''$ we have found a set that can play the r\^ole of $B$ in the Lemma's statement.

Finally, we need to consider the case where $T=+\infty$. One has $$ e^{-r\left(\upsilon\wedge n\right)}g(X_{\upsilon\wedge n})=\chi_{\left\{\upsilon<+\infty\right\}}e^{-r\left(\upsilon\wedge n\right)}g\left(X_{\upsilon\wedge n}\right)\longrightarrow e^{-r\upsilon}g(X_{\upsilon})\text{ as } n\rightarrow \infty \text{ $\PP$-a.s.}$$ for all $\cI$-valued $\upsilon$ that are almost surely finite (for, due to $X_{+\infty}=\Delta$ and $g\left(X_{+\infty}\right)$ by definition of the {\em cemetery} $\Delta$, one has $e^{-r\upsilon}g(X_{\upsilon})=\chi_{\left\{\upsilon <+ \infty\right\}}e^{-r\upsilon}g(X_{\upsilon})$). 
By our assumption of uniform integrability, we even have $L^1$-convergence in the previous convergence assertion (cf. e.g. Bauer \cite[Satz 21.4]{BauerMI}) and therefore obtain $$\lim_{n\rightarrow\infty}\EE\left[e^{-r\left(\upsilon\wedge n \right)}g(X_{\upsilon\wedge n})\right] = \EE\left[e^{-r \upsilon}g\left(X_{\upsilon}\right)\right].$$
Combining this result with the Lemma's statement for $T<+\infty$ (which has been proven before) we get the Lemma's estimate for $T=+\infty$, too.

\end{proof}

\begin{cor}[Formula for an option price using hitting times] \label{exerciseboundaryCorollary} Let $X$ be a $d$-dimensional basket with an associated risk-neutral measure $\PP$ and discount rate $r>0$. Consider a discrete subset $\cI\subset[0,+\infty)$. Suppose $g=(K-f)\vee 0$, and assume that the process $e^{-r\cdot}f(X_{\cdot})$ is a $\PP$-submartingale. Then one has \begin{eqnarray*}&&\sup_{\cG\subset \RR^d\times [0,T] \text{ measurable}}\EE \left[e^{-r\tau_\cG^\cI}g\left(X_{\tau_\cG^\cI}\right)\right] \\&=&\sup_{\tau \ \cI \cap [0,T]\text{-valued stopping time}}\EE \left[e^{-r\tau}g\left(X_\tau\right)\right]\end{eqnarray*} for all $T<+\infty$. If the the set $\left\{ e^{-r\tau}g\left(X_\tau\right) \ : \ {\tau\text{ $\cI$-valued stopping time} }\right\}$ is uniformly $\PP$-integrable, then the equation \begin{eqnarray*}&& \sup_{\cG\subset \RR^d\times [0,+\infty) \text{ measurable}}\EE \left[e^{-r\tau_\cG^\cI}g\left(X_{\tau_\cG^\cI}\right)\right] \\&=& \sup_{\tau \ \cI\text{-valued stopping time}}\EE \left[e^{-r\tau}g\left(X_\tau\right)\right]\end{eqnarray*} holds.

\end{cor}

Without going into detail we remark that similar results can be obtained for non-discrete $\cI$ as well, as proven by N El Karoui \cite{ElK}.

\begin{Def} Let $\cI\subset[0,+\infty)$ be discrete, $\cG\subset \RR^d\times [0,+\infty)$ and $G\subset \RR^d$ measurable, and $X$ a $d$-dimensional Markov basket with an associated family of risk-neutral measures $\PP^\cdot$ and discount rate $r>0$. We define 
$$V_{\cG,X}^\cI:(x,t)\mapsto\left\{\begin{array}{*{2}{c}} e^{rt}\EE^{(x,t)}\left[e^{-r\tau_\cG^\cI}g\left(X_{\tau_\cG^\cI}\right)\right],&(x,t)\notin \cG, \\ g(x), & (x,t)\in \cG.\end{array}\right.$$ as well as $$V_{G,X}^\cI:x\mapsto\left\{\begin{array}{*{2}{c}}\EE^x\left[e^{-r\tau_G^\cI}g\left(X_{\tau_G^\cI}\right)\right], &x\notin G, \\ g(x), & x\in G.\end{array}\right.$$ 

If $0\in\cI$, then we can simply write $$V_{\cG,X}^\cI:(x,t)\mapsto e^{rt}\EE^{(x,t)}\left[e^{-r\tau_\cG^\cI}g\left(X_{\tau_\cG^\cI}\right)\right]$$ and $$V_{G,X}^\cI:x\mapsto\EE^x\left[e^{-r\tau_G^\cI}g\left(X_{\tau_G^\cI}\right)\right].$$ 

Instead of $V_{\cG,X}^\cI(x,0)$, we shall often simply write $V_{\cG,X}^\cI(x)$. Also, the subscript $X$ will be dropped when no ambiguity can arise. Furthermore, $V_\cG^h$ and $V_G^h$ will be shorthand for $V_\cG^{h\NN_0}$ and $V_G^{h\NN_0}$, respectively.

\end{Def}

As another notational convention, let us from now on use $\sup_{\cG\subset \RR^d\times [0,T] }$ and $\sup_{G\subset \RR^d}$ to denote $\sup_{\cG\subset \RR^d\times [0,T] \text{ measurable}}$ and $\sup_{G\subset \RR^d\text{ measurable}}$, respectively.

The following Theorem is also a classical result from the theory of optimal stopping and Snell envelopes (cf. e.g. Neveu \cite[Proposition VI-2-8]{Nev} and El Karoui \cite[Th\'eor\`eme 2.31]{ElK})

\begin{Th}[Optimality of the immediate exercise region]\label{exerciseboundary} Let $X$ be a $d$-dimensional Feller basket with $\PP^\cdot$ being an associated family of risk-neutral measures and $r>0$ being the discount rate belonging to $\PP^\cdot$. Suppose $g=(K-f)\vee 0$, $\cI\subset[0,+\infty)$ is discrete, and $T\in[0,+\infty]$. Assume, moreover, that $e^{-r\cdot}f(X_{\cdot})$ is a $\PP^x$-submartingale for all $x\in\RR^d$. Define $$F^{\cI,T} =\left\{(x,t)\in\RR^d\times[0,T] \ : \sup_{\cG\subset \RR^d\times[0,T]}V_\cG^\cI(x,t) \leq g(x)\right\} $$ if $T<+\infty$ (we may drop the superscript $T$ wherever this is unambiguous) and else $$F^{\cI,+\infty} =\left\{(x,t)\in\RR^d\times[0,\infty) \ : \sup_{\cG\subset \RR^d\times[0,+\infty)}V_\cG^\cI(x,t) \leq g(x)\right\} $$ Then $$\forall x\in\RR^d \quad V_{F^{\cI,T}}^\cI(x,0)=U^{\cI}(T)(x)$$ if $T<+\infty$, and $V_{F^{\cI,+\infty}}^\cI(x,0)=U^{\cI}(+\infty)(x)$ for all $x\in\RR^d$ such that the set $\left\{ e^{-r\tau}g\left(X_\tau\right) \ : \ {\tau\text{ $\cI$-valued stopping time} }\right\}$ is uniformly $\PP^x$-integrable.
\end{Th}
\begin{proof} Let $T<+\infty$. Using Corollary \ref{exerciseboundaryCorollary} and recalling the definition of $U^{s\NN_0}$, all we have to show is $$\forall x\in\RR^d \quad V_{F^{\cI,T}}^\cI(x,0)= \sup_{\cG\subset \RR^d\times [0,T] }\EE^{(x,0)} \left[e^{-r\bar\tau_\cG^\cI}g\left(\bar\tau_\cG^\cI\right)\right]$$ (where we recall that $\bar\tau_\cG^\cI\leq \tau_\cG^\cI$ denotes the first nonnegative entry time into $\cG$). However, after exploiting the special particular shape of $F^\cI$, we can -- due to the boundedness of $g\leq K$ which yields $V_\cG^\cI\leq K$ for all $\cG$ which allows us to apply Lebesgue's Dominated Convergence Theorem -- swap $\sup$ and $\EE$ to get for all $x\in\RR^d$, \begin{eqnarray*}V_{F^\cI}^\cI(x,0) &=& \EE^{(x,0)}\left[e^{-r\bar\tau_{F^\cI}^\cI} g\left(X_{\bar\tau_{F^\cI}^\cI}\right) \right] \\ &=& \EE^{(x,0)}\left[e^{-r\bar\tau_{F^\cI}^\cI} \sup_{\cG\subset \RR^d\times[0,T]}V_\cG^\cI\left(X_{\bar\tau_{F^\cI}^\cI},{\bar\tau_{F^\cI}^\cI}\right) \right] \\ &=& \EE^{(x,0)}\left[e^{-r\bar\tau_{F^\cI}^\cI} \sup_{\cG\subset \RR^d\times[0,T]} e^{r\bar\tau_{F^\cI}^\cI}\EE^{\left(X_{\bar\tau_{F^\cI}^\cI},{\bar\tau_{F^\cI}^\cI}\right)} e^{-r\bar\tau_{\cG}^\cI} g\left(X_{\bar\tau_{\cG}^\cI}\right)\right] \\ &=& \sup_{\cG\subset \RR^d\times [0,T]} \EE^{(x,0)}\left[\EE^{\left(X_{\bar\tau_{F^\cI}^\cI},{\bar\tau_{F^\cI}^\cI}\right)}e^{-r\bar\tau_{\cG}^\cI}g\left(X_{\bar\tau_\cG^\cI}\right) \right]\end{eqnarray*} (where for notational convenience $\bar\tau_\cG^\cI$ should denote the first nonnegative entry time into $\cG$). Now, let us use the strong Markov property of the Feller process $X$, and for this purpose, let $\bar\theta$ denote the shift operator on the space-time path space $D\left([0,+\infty),\RR^d\times [0,+\infty)\right)$ (which is the set of all c\`adl\`ag functions from $[0,+\infty)$ into $\RR^d\times [0,+\infty)$ -- recall that all Feller processes have a c\`adl\`ag modification). We obtain \begin{eqnarray*}V_{F^\cI}^\cI(x,0) &=& \sup_{\cG\subset \RR^d\times [0,T]} \EE^{(x,0)}\left[\EE^{\left(X_{\bar\tau_{F^\cI}^\cI},{\bar\tau_{F^\cI}^\cI}\right)}e^{-r\bar\tau_{\cG}^\cI}g\left(X_{\bar\tau_\cG^\cI}\right) \right]\\ &=& \sup_{\cG\subset \RR^d\times [0,T]} \EE^{(x,0)}\left[\EE^{\left(X_{0},{0}\right)}\left[\left.e^{-r\bar\tau_{\cG}^\cI\circ \bar\theta_{\bar\tau_{F^\cI}^\cI}}g\left(X_{\bar\tau_\cG^\cI}\circ \bar\theta_{\bar\tau_{F^\cI}^\cI}\right)\right|\cF_{\bar\tau_{F^\cI}^\cI}\right] \right] \\ &=& \sup_{\cG\subset \RR^d\times [0,T]} \EE^{(x,0)}\left[e^{-r\bar\tau_\cG^\cI\circ \bar\theta_{\bar\tau_{F^\cI}^\cI}}g\left(X_{\bar\tau_\cG^\cI\circ \bar\theta_{\bar\tau_{F^\cI}^\cI}}\right) \right]\end{eqnarray*} 

But $e^{-r\cdot}g(X_{\cdot})$ is a $\PP^x$-supermartingale for all $x\in\RR^d$, therefore by Doob's Optional Stopping Theorem, $\left(e^{-r\upsilon}g\left(X_{\upsilon}\right)\right) _{\upsilon\in\left\{\bar\tau_\cG^\cI\wedge T',\bar\tau_\cG^\cI\circ \bar\theta{\tau_{F^\cI}^\cI}\wedge T'\right\}}$ must also be a $\PP^x$-submartingale for all $x\in\RR^d$ and $T'\in(0,\infty)$ (note that $\bar\tau_\cG^\cI\leq \bar\tau_\cG^\cI\circ \bar\theta{\tau_{F^\cI}^\cI}$ a.s. because of the fact that $\bar\theta$ is the shift operator for the space-time process $(t,X_t)_{t\geq 0}$, rather than simply for $X$). Letting $T'$ tend to infinity, we can employ Lebesgue's Dominated Convergence Theorem (as $g\leq K$ yields $e^{-r\upsilon}g\left(X_{\upsilon}\right)\leq K\in L^1(\PP^x)$ for $\upsilon\in\left\{\bar\tau_\cG^\cI\wedge T',\bar\tau_\cG^\cI\circ \bar\theta_{\bar\tau_{F^\cI}^\cI}\wedge T'\right\}$ for all $T'\in(0,+\infty)$ and $x\in\RR^d$) in order to get that the expected value of $e^{-r\bar\tau_\cG^\cI\circ \bar\theta_{\bar\tau_{F^\cI}^\cI}}g\left(X_{\bar\tau_\cG^\cI\circ \bar\theta_{\bar\tau_{F^\cI}^\cI}}\right)$ is always greater than or equal to the expectation of $e^{-r\bar\tau_\cG^\cI}g\left(X_{\bar\tau_\cG^\cI}\right)$.
Hence \begin{eqnarray*}V_{F^\cI}^\cI(x,0) &=& \sup_{\cG\subset \RR^d\times [0,T]} \EE^{(x,0)}\left[e^{-r\bar\tau_\cG^\cI\circ \bar\theta_{\bar\tau_{F^\cI}^\cI}}g\left(X_{\bar\tau_\cG^\cI\circ \bar\theta_{\bar\tau_{F^\cI}^\cI}}\right) \right] \\ &\geq & \sup_{\cG\subset \RR^d\times [0,T]} \EE^{(x,0)}\left[e^{-r\bar\tau_\cG^\cI}g\left(X_{\bar\tau_\cG^\cI}\right) \right]\\ &=& \sup_{\cG\subset \RR^d\times [0,T]} \EE^{(x,0)}V_\cG^\cI(x,0).\end{eqnarray*} 

The case $T=+\infty$ can be dealt with analogously.

\end{proof}

\begin{rem}\label{American_region} Whenever Lemma \ref{dyadic limits generic} may be applied, the immediate exercise region of an American option is the intersection of the immediate exercise regions of the corresponding approximating Bermudan options: \begin{eqnarray*} && \left\{(x,t)\in\RR^d\times[0,T] \ : \sup_{\cG\subset \RR^d\times[0,T]}V_\cG^{[0,+\infty]}(x,t) \leq g(x)\right\}\\ &=& \bigcap_{s>0}F^{s\NN_0,T} =\bigcap_{m\in\NN}F^{2^{-m}\NN,T}, \end{eqnarray*} and its optimality follows from the general theory of optimal stopping (cf. Neveu \cite{Nev}, Griffeath and Snell \cite{Snell}).
\end{rem}

\begin{lem}[Time-stationarity of immediate exercise regions for perpetual Bermudans]\label{time-stationarity} Let $X$ be a L\'evy basket with $\PP^\cdot$ being an associated family of probability measures and discount rate $r>0$. Then for all $s>0$ we have $$U^{s\NN_0}(+\infty)(x)=V_{\left\{x\in\RR^d\ : \ U^{s\NN_0}(+\infty)(x)\leq g(x)\right\}}^{s\NN_0}(x)$$ for all $x\in\RR^d$ satisfying the condition that $\left\{e^{-r\tau}g\left(X_\tau\right) \ : \ \tau\text{ $s\NN_0$-valued stopping time}\right\}$ is uniformly $\PP^x$-integrable.
\end{lem}
\begin{proof} Consider an integer $n\in \NN_0$, and an $x\in\RR^d$ such that $\sup\left\{e^{-r\tau}g\left(X_\tau\right) \ : \ \tau\text{ stopping time}\right\}$ is $\PP^x$-integrable. Then we shift the time scale by $ns$ to get\begin{eqnarray*}&& e^{rns}\EE^{(x,ns)} \left[\sup_{\cG\subset \RR^d\times [0,+\infty)}e^{-r\tau_\cG^{s}} g\left(X_{\tau_\cG^{s}}\right)\right] \\&=& e^{rns}\EE^{(x,ns)} \left[\sup_{\cG\subset \RR^d\times [ns,+\infty)}e^{-r\tau_\cG^{s}} g\left(X_{\tau_\cG^{s}}\right)\right] \\&=& e^{rns}\EE^{(x,ns)} \left[\sup_{\cG'\subset \RR^d\times [0,+\infty)}e^{-r\left(\tau_{\cG'}^{s}\circ\theta_{ns}+ns\right)} g\left(X_{\tau_{\cG'}^{s}}{\circ\theta_{ns}}\right)\right] \\&=& \EE^{(x,0)} \left[\sup_{\cG'\subset \RR^d\times [0,+\infty)}e^{-r\tau_{\cG'}^{s}} g\left(X_{\tau_{\cG'}^{s}}\right)\right] 
\end{eqnarray*} where $\theta$ denotes the shift operator on the space (as opposed to space-time) path space $D\left([0,+\infty),\RR^d\right)$. Because of the boundedness of $g\leq K$ which entitles us to apply Lebegue's Dominated Convergence Theorem, we may swap $\sup$ and $\EE$ to obtain \begin{eqnarray*} && e^{rns}\EE^{(x,ns)} \left[\sup_{\cG\subset \RR^d\times [0,+\infty)}e^{-r\tau_\cG^{s}} g\left(X_{\tau_\cG^{s}}\right)\right]\\ &=& \sup_{\cG\subset \RR^d\times [0,+\infty)}e^{rns}\EE^{(x,ns)} \left[e^{-r\tau_\cG^{s}} g\left(X_{\tau_\cG^{s}}\right)\right]\\ &=& \sup_{\cG\subset \RR^d\times [0,+\infty)}V_\cG^{s\NN_0}(x,ns)\end{eqnarray*} for all $n\in\NN_0$. Thus we conclude $$\sup_{\cG\subset \RR^d\times [0,+\infty)}V_\cG^{s\NN_0}(x,t)=\sup_{\cG\subset \RR^d\times [0,+\infty)}V_\cG^{s\NN_0}(x,0)$$ for all $t\in s\NN_0$. If we insert this equality fact into the definition of $F^{s\NN_0}$, we see that the condition determining whether a pair $(x,t)$ belongs to $F^{s\NN_0}$ does not depend on $t$. On the other hand, by Corollary \ref{exerciseboundaryCorollary}, $$\sup_{\cG\subset \RR^d\times [0,+\infty)}V_\cG^{s\NN_0}(x,0)=U^{s\NN_0}(+\infty)(x),$$ and the left hand side equals -- by our previous observations in this proof -- the term featuring in the definition of $F^{s\NN_0}$.

\end{proof}
Summarising the previous deliberations, we deduce that under the assumptions of the previous Lemmas, the expected payoff of a perpetual Bermudan option of mesh size $s>0$ equals $$U^{s\NN_0}(+\infty)(\cdot)= \EE^\cdot \left[e^{-r\tau_{G^s}^ s}g\left(X_{\tau_{G^s}^s}\right)\right] ,$$ where $G^s:=G^s_g:=\left\{x\in\RR^d\ : \ U^{s\NN_0}(+\infty)(x)\leq g(x)\right\}= \left\{x\in\RR^d\ : \ \sup_{H\subseteq\RR^d\atop \text{measurable}} \EE^{x}\left[e^{-r\tau_H^s}g\left(X_{\tau_H^s}\right)\right]\leq g(x)\right\}$.

\begin{rem} As we have already observed in Remark \ref{American_region}, under the assumptions of Lemma \ref{dyadic limits generic} and Theorem \ref{exerciseboundary}, the region $G^\ast_g:=\bigcap_{n\in\NN}G^{2^{-n}}_g= \bigcap_{s>0}G^{s}_g$ (where $G^{s}$ denotes the immediate exercise region of the perpetual Bermudan option with payoff function $g$ and exercise mesh $s$) will be the optimal exercise region of the American perpetual option with payoff function $g$. 

Furthermore, whenever $G\subseteq G^\ast_g$, then the immediate exercise region for the perpetual Bermudan option with payoff function $g\chi_{G}$ and exercise mesh $s>0$, denoted $G^{s}_{\chi_Gg}$, will again be $G$ (and will thus no longer depend on $s$): This follows, due to $G^\ast_g=\bigcap_{n\in\NN}G^{2^{-n}}_g=\bigcap_{s>0}G^s_g$, from the more general fact that whenever $G\subseteq G^s_g$ for some $s>0$, then $G^s_{\chi_Gg}=G$. (To prove this fact, simply observe that on the one hand $G^s_{\chi_Gg}\supseteq G^s_g$, since $\sup_{H\subseteq\RR^d\atop \text{measurable}} \EE^{\cdot}\left[e^{-r\tau_H^s}(\chi_{G}g)\left(X_{\tau_H^s}\right)\right]\leq \sup_{H\subseteq\RR^d\atop \text{measurable}} \EE^{\cdot}\left[e^{-r\tau_H^s}g\left(X_{\tau_H^s}\right)\right]\leq g(x)$. Hence, $G^s_{\chi_Gg}\supseteq G$, since $G^s_g\supseteq G$ by assumption. On the other hand, the immediate exercise region $G^s_{\chi_Gg}$ must be inside $G$, for the payoff function $\chi_Gg$ trivially vanishes outside $G$.) 

\end{rem}

\begin{lem} \label{convolveproject}Let us fix a L\'evy basket with an associated family of risk-neutral probability measures $\PP^\cdot$ and discount rate $r>0$, as well as a region $G\subset \RR^d$ and a real number $s>0$. Then we have $$\forall x\not\in G \quad V_G^s(x)=e^{-rs}\PP_{X_0-X_s}\ast V_G^s(x).$$ In particular, using Lemma \ref{time-stationarity}, one has the following equation for the expected perpetual Bermudan option payoff: $$\forall x\not\in G \quad U^{s\NN_0}(+\infty)(x)=e^{-rs}\PP_{X_0-X_s}\ast U^{s\NN_0}(+\infty)(x).$$
\end{lem}
\begin{proof} Using the Markov property of $X$, denoting by $\theta$ the shift operator on the path space $D\left([0,+\infty),\RR^d\right)$ of a L\'evy process $X$, and taking into account the fact that $\tau_G^s>0$ (ie $\tau_G^s\geq s$) in case $x\not\in G$, we obtain:
\begin{eqnarray*}\forall x\not\in G \quad V_G^s(x) &=& e^{-rs}\EE^x \EE^{x}\left[ e^{-r\left(\tau_G^s-s\right)}g\left(X_{\tau_G^s}\right) \cdot \chi_{\{\tau_G^s\geq s\}}\right] \\ &=& e^{-rs}\EE^x \EE^{x}\left[ e^{-r\tau_G^s\circ \theta_s}g\left(X_{\tau_G^s}\circ\theta_s \right) \cdot \chi_{\{\tau_G^s\geq s\}}\right] \\ &=& e^{-rs} \EE ^x\EE^{x}\left[ e^{-r\tau_G^s\circ \theta_s}g\left(X_{\tau_G^s}\circ\theta_s \right) \left( \chi_{\{\tau_G^s\geq s\}} + \chi_{\{\tau_G^s<s\}}\right)\right] \\ &=& e^{-rs} \EE^x \EE^x \left[ e^{-r\tau_G^s \circ \theta_s}g\left(X_{\tau_G^s}\circ\theta_s\right) | \cF_{s} \right]\\ &=& e^{-rs} \EE^x \EE^{X_s} e^{-r\tau_G^s}g\left(X_{\tau_G^s}\right) = e^{-rs}\EE ^xV_G^s\left(X_s\right) \\ &=& e^{-rs} \PP_{x-X_s}\ast V_G^s(x) \\ &=& e^{-rs} \PP_{X_0-X_s}\ast V_G^s(x). \end{eqnarray*} 
\end{proof}


\section{Natural scaling for continuity corrections of American perpetuals}

Adopting the terminology of Broadie, Glasserman and Kou \cite{BGK}, we shall refer to the difference between an American and the corresponding Bermudan options (on the same basket and with the same payoff function) as a ``continuity correction''.

\subsection{Continuity corrections}

Let $X$ be a $d$-dimensional Markov basket with an associated family of risk-neutral measures $\PP^\cdot$ and discount rate $r>0$. Consider a (non-dividend paying) perpetual Bermudan options with exercise mesh size $s>0$, a measurable set $G\subset\RR^d$, and a bounded nonnegative measurable function $g:\RR^d\rightarrow \RR_{\geq 0}$ (the assumption of boundedness applies e.g. to the case of put options). The expected payoff of this option with respect to the exercise region $G$ for logarithmic start prices $x\in \complement G$ is then given by $$V_{G}^s(x) = \EE^x\left[e^{-r\tau_{G}^{s}}g\left(X_{\tau_{G}^{s}}\right)\right].$$

\begin{rem} We shall consider the problem of pricing an American perpetual option with payoff function $g$. Let $G$ be its optimal exercise region. Then it is enough to compute the price for an American perpetual option with payoff function $g\cdot \chi_G$. As was shown in Lemma \ref{dyadic limits generic}, if $X$ has a modification with continuous paths and $g$ is continuous, then the American price will be the limit, as $s\downarrow 0$, of the expected payoff of a Bermudan option with payoff function $g\cdot \chi_G$ and exercise mesh size, viz. $V_G^s$. Since $G$ was the optimal exercise region of the American option, $G$ will be, independent of $s>0$, the optimal exercise region of a Bermudan option with payoff function $g\cdot \chi_G$.
\end{rem}

We start with the definition of what will be used as a measure for continuity corrections.

\begin{Def} Let $X$ be a $d$-dimensional Markov basket with an associated family of risk-neutral measures $\PP^\cdot$ and discount rate $r>0$. Consider a bounded continuous payoff function $g:\RR^d\rightarrow \RR_{\geq 0}$. We define
$$ \forall s>0 \forall x\in\RR^d\quad \rho(s)(x):= \lim_{t\downarrow 0}\EE^x\left[e^{-r\tau_{G}^{t\NN}} - e^{-r\tau_{G}^{s\NN}}\right] =\left. \EE^x\left[ e^{-r\tau_{G}^{t\NN}}\right] \right|^{t\downarrow 0}_{t=s}$$ whenever this limit exists.
\end{Def}

\begin{rem} 
\begin{enumerate}
\item If $X$ has a modification with continuous paths, then (readily due to Lemma \ref{dyadic limits generic} and monotone convergence), $$ \rho(s)=\EE^\cdot\left[\sup_{m\in\NN}e^{-r\tau_G^{s\cdot 2^{-m}}}- e^{-r\tau_{G}^t}\right] = \sup_{m\in\NN} \EE^\cdot\left[e^{-r\tau_G^{s\cdot 2^{-m}}}- e^{-r\tau_{G}^t}\right] ,$$ hence the limit in the definition of $\rho(s)$ exists, whence $\rho(s)$ is well-defined.
\item $\rho(s)$ provides an estimate for the continuity correction. Since $g$ is bounded and nonnegative, $$ \forall x\in\complement G \quad \left| V^t_G(x)- V^s_G(x)\right|\leq \sup_G g \cdot \EE^x\left[e^{-r\tau_{G}^t}- e^{-r\tau_{G}^s}\right],$$ so \begin{equation}\label{cc vindic} \left| \lim_{t\downarrow 0}V^t_G- V^s_G\right|\leq \sup_G g\cdot \rho(s) \text{ on }\complement G.\end{equation}
\item $\rho(s)$, for any $s>0$, depends only on $G$, not on the payoff function $g$ itself.
\end{enumerate}
\end{rem}

\subsection{Formulae for $\rho(s)$ in dimension one}

Formulae for continuity corrections for barrier options in the one-dimensional Black-Scholes model have already been derived by other authors (e.g. Broadie, Glassermann and Kou \cite{BGKbarrier}, H\"orfelt \cite{Hoe}, Howison \cite{How}, as well as Howison and Steinberg \cite{HS}). We will nevertheless, for the sake of illustration, show how our approach applies to this (compared to other examples, of course, very simple) setting. This should be seen as a motivation for the proof of Theorem \ref{phase} and its Corollaries, where higher dimensions are studied.

\begin{Th}\label{phase1d} Consider a $1$-dimensional L\'evy basket $X$ with an associated family of risk-neutral measures $\PP^\cdot$ and discount rate $r>0$. 
Let $G=(-\infty,\gamma)$ for some $\gamma\in\RR$. Then one has for all $s>0$ the relations \begin{eqnarray} \rho(s)(\gamma)\nonumber &=& \left.\sum_{n=1}^\infty e^{-rnt} \PP^{\gamma} \left[ \bigcap_{i=1}^{n-1}\left\{ X_{it}\geq \gamma \right\}\cap\left\{X_{nt} < \gamma \right\}\right]\right|^{t\downarrow 0}_{t=s}\nonumber \\ &=& \left[\exp\left(-\sum_{n=1}^\infty\frac{e^{-rnt}}{n} \PP^0\left\{X_{nt}< 0\right\}\right) \right]^{t=s}_{t\downarrow 0}\label{phase1d identity}\end{eqnarray} (and the limit on the right exists whenever $\rho(s)(\gamma)$ is well-defined).
\end{Th}
\begin{cor} \label{BlackScholes1d} Let $G=(-\infty,\gamma]$ or $G=(-\infty,\gamma)$ for some $\gamma\in\RR$ and assume $(X_t)_{t\geq 0}=\left(X_0+\sigma\cdot B_t+\left(r-\frac{\sigma^2}{2}\right)\right)_{t\geq 0}$, in words: $X$ is the logarithmic price process of the one-dimensional Black-Scholes model with constant volatility $\sigma$ and discount rate $r>0$. Then, whenever $\mu:=r-\frac{\sigma^2}{2}\geq 0$, there exist constants $c_0, C_0>0$ such that for all sufficiently small $s>0$, \begin{eqnarray*} c_0s^\frac{1}{2}\leq \rho(s)(\gamma) \leq C_0s^\frac{1}{2\sqrt{2}}\end{eqnarray*} If both $\mu\leq 0$ and $r> \frac{\mu^2}{2\sigma^2}$, there exist constants $c_1, C_1>0$ such that for all sufficiently small $s>0$, \begin{eqnarray*}c_1s^\frac{1}{\sqrt{2}} \leq \rho(s)(\gamma)\leq C_1s^\frac{1}{{2}}.\end{eqnarray*} 
\end{cor}

\begin{rem} Although computing the constants $c_0,C_0,c_1,C_1$ explicitly is possible, we refrain from doing so for the moment, as it is not required to find the right scaling for an extrapolation for $V_G^s$ from $s>0$ to $s=0$ and it would not provide any additional useful information for our extrapolation purposes. The same remark applies to all examples and generalisations that are studied subsequently.
\end{rem}

\begin{proof}[Proof of Theorem \ref{phase1d}] The Theorem follows from a result by Feller \cite[p. 606, Lemma 3]{F} on processes with stationary and independent increments. For, if we define $$\forall s>0\quad\forall q\in[0,1) \quad \xi(q,s):= \sum_{n=1}^\infty q^n \PP^{0} \left[ \bigcap_{i=1}^{n-1}\left\{ X_{is}\geq 0\right\}\cap\left\{X_{ns}< 0\right\}\right] ,$$ then Feller's identity \cite[p. 606, Lemma 3]{F} reads \begin{equation}\label{Feller} \forall s>0\quad\forall q\in[0,1)\quad -\ln \left(1-\xi(q,s)\right) = \sum_{n=1}^\infty\frac{q^n}{n} \PP^0\left\{X_{ns}< 0\right\} \end{equation} and holds whenever $X$ has stationary and independent increments, in particular for all L\'evy processes (note that our definition of a L\'evy process requires them to be Feller processes in addition). This entails \begin{eqnarray} \label{rhosimVG}\xi\left(e^{-rs},s\right)&=&\sum_{n=1}^\infty e^{-rns} \PP^{0} \left[ \bigcap_{i=1}^{n-1}\left\{ X_{is}\geq 0\right\}\cap\left\{X_{ns}< 0\right\}\right] \\&=& 1-\exp\left(-\sum_{n=1}^\infty\frac{e^{-rns}}{n} \PP^0\left\{X_{ns}< 0\right\}\right) ,\end{eqnarray} which is enough to prove the second identity (\ref{phase1d identity}) in the Theorem.
\end{proof}

\begin{proof}[Proof of Corollary \ref{BlackScholes1d}]
First, it makes no difference whether we consider $G_1=(-\infty,\gamma)$ or $G_1=(-\infty,\gamma]$, as $G_1\setminus G_0=\{\gamma\}$ which has capacity zero. In light of Theorem \ref{phase1d}, we shall show that if $\mu\geq 0$, there exist constants $c_0, C_0>0$ such that for all sufficiently small $s>0$, \begin{eqnarray*} c_0s^\frac{1}{2}\leq \exp\left(-\sum_{n=1}^\infty \frac{e^{-rns}}{n}\PP^{0} \left\{X_{ns}\leq 0\right\} \right) \leq C_0s^\frac{1}{2\sqrt{2}},\end{eqnarray*} and if both $\mu\leq 0$ and $r> \frac{\mu^2}{2\sigma^2}$, there exist constants $c_1, C_1>0$ such that for all sufficiently small $s>0$, \begin{eqnarray*}c_1s^\frac{1}{\sqrt{2}} \leq \exp\left(-\sum_{n=1}^\infty \frac{e^{-rns}}{n}\PP^{0} \left\{X_{ns}\leq 0\right\} \right) \leq C_1s^\frac{1}{{2}}.\end{eqnarray*} Now, the scaling invariance of Brownian motion yields for all $n\in\NN$ and $s>0$: \begin{eqnarray} \label{scaling}\PP^0\left\{X_{ns}\leq 0\right\}&=& \PP^0\left\{B_{ns}\leq -\frac{\mu}{\sigma} ns\right\}=\PP^0\left\{ B_1\leq -\frac{\mu}{\sigma}(ns)^{1/2}\right\} \nonumber \\ &=& (2\pi)^{-1/2}\int_{-\infty}^{-\frac{\mu}{\sigma}(ns)^{1/2}} \exp\left(\frac{-x^2}{2}\right)dx .\end{eqnarray} We divide the remainder of the proof, which will essentially consist in finding estimates for the right hand side of the last equation, into two parts according to the sign of $\mu$.\\
{\em Case I: $\mu\geq 0$.} In this case we use the estimates $$\forall x\leq 0\forall y\leq 0\quad -y^2-x^2\leq -\frac{|x+y|^2}{2} \leq -\frac{y^2}{2}-\frac{x^2}{2}, $$ thus \begin{eqnarray*}&&\forall y\leq 0 \\ &&e^{-y^2}\int_{-\infty}^0 e^{-x^2}dx \leq\int_{-\infty}^0 \exp\left(-\frac{|x+y|^2}{2}\right) dx\leq e^{-y^2/2}\int_{-\infty}^0 e^{-x^2/2}dx ,\end{eqnarray*} hence by transformation for all $y\leq 0$
\begin{equation}\label{elestim}\frac{ \sqrt{\pi}}{2} e^{-y^2}\leq \int_{-\infty}^y \exp\left(-\frac{x^2}{2}\right) dx \leq \sqrt{\frac{\pi}{2}}e^{-\frac{y^2}{2}} .\end{equation} Due to equation (\ref{scaling}), this entails for all $n\in\NN$, $s>0$, $\mu\geq0$ (if we insert $-\frac{\mu}{\sigma}(ns)^{1/2}$ for $y$) $$\frac{e^{-\frac{\mu^2}{\sigma^2} ns}}{2\sqrt{2}}\leq \PP^0\left\{ X_{ns}\leq 0 \right\}\leq \frac{e^{-\frac{\mu^2 ns }{2\sigma^2}}}{2}.$$ Therefore for arbitrary $r,s>0$, \begin{equation}\label{geoharmonic}\frac{1}{2\sqrt{2}}\sum_{n=1}^\infty \frac{e^{-ns\left(r+\frac{\mu^2}{\sigma^2}\right)}}{n} \leq \sum_{n=1}^\infty \frac{e^{-rns}}{n}\PP^0\left\{X_{ns}\leq 0\right\} \leq \frac{1}{2}\sum_{n=1}^\infty \frac{e^{-ns\left(r+\frac{\mu^2}{2\sigma^2}\right)}}{n} \end{equation} The sums in equation (\ref{geoharmonic}) have got the shape of $\sum q^n/n$ for $q<1$. 
Now one performs a standard elementary computation on this power series: \begin{eqnarray}\label{geoharmo1} \sum_{n=1}^\infty \frac{q^n}{n}&=&\sum_{n=0}^\infty \int_0^q r^n dr= \int_0^q\sum_{n=0}^\infty r^n dr = \int_0^q \frac{1}{1-r} dr=-\ln(1-q) ,\end{eqnarray} which immediately gives \begin{eqnarray*} &&\left(1-e^{-s\left(r+\frac{\mu^2}{2\sigma^2}\right)}\right)^{1/2} \\ &\leq &\exp\left(- \sum_{n=1}^\infty \frac{e^{-rns}}{n}\PP^0\left\{X_{ns}\leq 0\right\}\right)\\ & \leq &\left(1-e^{-s\left(r+\frac{\mu^2}{\sigma^2}\right)}\right)^\frac{1}{2\sqrt{2}} \end{eqnarray*} when applied to equation (\ref{geoharmonic}). Due to de l'Hospital's rule, the differences in the brackets on the left and right hand sides of the last estimate behave asymptotically like $s$ when $s\downarrow 0$. This is sufficient to prove the estimate in the Corollary for the case of $\mu\geq 0$. \\
{\em Case II: $\mu\leq 0$ and $r>\frac{\mu^2}{ 2\sigma^2 }$.} In that case we employ the estimates $$\forall x\leq 0\quad \forall y\leq 0\quad -\frac{x^2}{2}-\frac{y^2}{2}\leq -\frac{|x-y|^2}{2} \leq -\frac{x^2}{4}+\frac{y^2}{2} $$ and proceed analogously to Case I, to obtain \begin{equation}\label{elestimII}\sqrt{\frac{ {\pi}}{2}} e^\frac{-y^2}{2}\leq \int_{-\infty}^{-y} \exp\left(-\frac{x^2}{2}\right) dx \leq \sqrt{{\pi}}e^{\frac{y^2}{2}} .\end{equation} In the special case of $y:= \frac{\mu}{\sigma}(ns)^{1/2}\leq 0$, this leads to the estimate in the statement of the Corollary via $$ \frac{ e^{-\frac{\mu^2}{2\sigma^2}ns} }{2} \leq \PP^0\left\{X_{ns}\leq 0\right\}\leq \frac{ e^{\frac{\mu^2}{2\sigma^2}ns} }{\sqrt{2}}$$ and \begin{equation}\label{}\frac{1}{2}\sum_{n=1}^\infty \frac{e^{-ns\left(r+\frac{\mu^2}{2\sigma^2}\right)}}{n} \leq \sum_{n=1}^\infty \frac{e^{-rns}}{n}\PP^0\left\{X_{ns}\leq 0\right\} \leq \frac{1}{\sqrt{2}}\sum_{n=1}^\infty \frac{e^{-ns\left(r-\frac{\mu^2}{2\sigma^2}\right)}}{n} .\end{equation} \\
Therefore in case $\mu=0$ the scaling exponent is exactly $\frac{1}{2}$.
\end{proof}

The identity (\ref{phase1d identity}) of Theorem \ref{phase1d} can be used to derive estimates in the spirit of Corollary \ref{BlackScholes1d} in more general situations. We will illustrate this by means of the following example:

\begin{ex}[Merton's jump-diffusion model with positive jumps and ``moderate'' volatility] Suppose the logarithmic price process $X$ is governed by an equation of the form $$ \forall t\geq 0 \quad X_t=X_0+\alpha t+ \beta Z_t+ \sigma B_t $$ where $\alpha\in\RR$, $\beta,\sigma>0$, $Z$ is the Poisson process (thus, in this setting, only positive jumps are allowed for simplicity) for the parameter $1$ and $B$ a normalised one-dimensional Brownian motion, and the stochastic processes $B$ and $Z$ are assumed to be independent. Let $\PP^\cdot$ be an associated family of risk-neutral measures and $r>0$ the discount rate. In order to employ (\ref{phase1d identity}), we shall compute the sum $\sum_{n=0}^\infty\frac{e^{-rns}}{n}\PP^0\left\{X_{ns} < 0\right\}$ for all $s>0$. Since $\PP^0\left\{X_{ns}< 0\right\}=\PP^0\left\{\frac{X_{ns}}{\beta}< 0\right\}$ for arbitrary $n,s$ we may without loss of generality take $\beta =1$. Given $\alpha\leq r-\ln\EE^0Z_1$, the process $\exp\left(X_t-rt\right)_{t\geq 0}$ will be a martingale for the unique $\sigma$ satisfying \begin{align*}1=\EE^0\left[\exp\left(\alpha -r+\frac{\sigma^2}{2}+ {Z_1}\right)\right].\end{align*} Let us, in addition, assume $\alpha\geq 0$. Then $\exp\left(X_t-rt\right)_{t\geq 0}$ being a martingale implies $$\sigma=\sqrt{2\left(r-\alpha+ \ln \EE^0\left[e^{Z_1}\right]\right)}\leq \sqrt{2\left(r+ \ln \EE^0\left[e^{Z_1}\right]\right)}.$$ Now, by definition of the Poisson distribution together with the symmetry and scaling invariance of Brownian motion \begin{eqnarray} \label{jumpsum1} \nonumber &&\sum_{n=0}^\infty\frac{e^{-rns}}{n}\PP^0\left\{X_{ns} < 0\right\} \\ \nonumber &=& \sum_{n=0}^\infty \sum_{k=0}^\infty e^{-ns}\frac{(ns)^k}{k!}\cdot \frac{e^{-rns}}{n}\PP^0\left\{\sigma B_{ns}< -\alpha ns-k\right\} \\ &=& \sum_{n=0}^\infty \sum_{k=0}^\infty e^{-ns}\frac{(ns)^k}{k!}\cdot \frac{e^{-rns}}{n}\underbrace{\PP^0\left\{B_{1}< -\frac{\alpha}{\sigma} (ns)^{\frac{1}{2}}-\frac{k}{\sigma}(ns)^{-\frac{1}{2}}\right\}}_{=(2\pi)^{-1/2}\int_{-\infty}^{-\frac{\alpha}{\sigma} (ns)^{\frac{1}{2}}-\frac{k}{\sigma}(ns)^{-\frac{1}{2}}} \exp\left(\frac{-x^2}{2}\right)dx} \end{eqnarray} (with the convention that $0^0=1$). Now let us first of all try and find estimates for the probability in the last line. By equation (\ref{elestim}) applied to $y:=-\frac{\alpha}{\sigma} (ns)^{\frac{1}{2}}-\frac{k}{\sigma}(ns)^{-\frac{1}{2}}\leq 0$, \begin{eqnarray*} \frac{ e^{ -\left(\frac{\alpha}{\sigma} (ns)^{\frac{1}{2}}+\frac{k}{\sigma}(ns)^{-\frac{1}{2}} \right)^2 } }{2\sqrt{2}} &\leq& \PP^0\left\{B_{1}< -\frac{\alpha}{\sigma} (ns)^{\frac{1}{2}}-\frac{k}{\sigma}(ns)^{-\frac{1}{2}}\right\} \\ &\leq& \frac{ e^\frac{-\left(\frac{\alpha}{\sigma} (ns)^{\frac{1}{2}}+\frac{k}{\sigma}(ns)^{-\frac{1}{2}} \right)^2}{2} }{{2}} \end{eqnarray*} which yields, using the abbreviation $\alpha':=\frac{\alpha}{\sigma}$, \begin{eqnarray*} \frac{ e^{ -{\alpha'}^2 ns -2 \frac{\alpha'}{\sigma}k - \frac{k^2}{\sigma^2ns} } }{2\sqrt{2}} &\leq& \PP^0\left\{B_{1} < -\frac{\alpha}{\sigma} (ns)^{\frac{1}{2}}-\frac{k}{\sigma}(ns)^{-\frac{1}{2}}\right\} \\ &\leq& \frac{ e^{ -\frac{{\alpha'}^2}{2} ns -\frac{\alpha '}{\sigma}k -\frac{k^2}{2\sigma^2ns} } }{2} , \end{eqnarray*} so \begin{eqnarray} \label{lowerjump1} \frac{ e^{ -{\alpha'}^2 ns -k\left(2 \frac{\alpha'}{\sigma} + \frac{k}{\sigma^2ns}\right) } }{2\sqrt{2}} &\leq& \PP^0\left\{B_{1} <-\frac{\alpha}{\sigma} (ns)^{\frac{1}{2}}-\frac{k}{\sigma}(ns)^{-\frac{1}{2}}\right\} \\ \nonumber &\leq& \frac{ e^{ -\frac{{\alpha'}^2}{2} ns -\frac{\alpha '}{\sigma}k } }{2}. \end{eqnarray} Thus, we can perform the following estimates to derive an upper bound of the sum in (\ref{jumpsum1}): \begin{eqnarray*} && \sum_{n=0}^\infty \sum_{k=0}^\infty e^{-ns}\frac{(ns)^k}{k!}\cdot \frac{e^{-rns}}{n}\PP^0\left\{B_{1}< -\frac{\alpha}{\sigma} (ns)^{\frac{1}{2}}-\frac{k}{\sigma}(ns)^{-\frac{1}{2}}\right\} \\ &\leq& \frac{1}{2} \sum_{n=0}^\infty \frac{e^{-ns\left( 1+r+\frac{{\alpha'}^2}{2}\right)}}{n}\sum_{k=0}^\infty\frac{1}{k!}\left(e^\frac{-\alpha '}{\sigma}\cdot ns\right)^k \\ &=& \frac{1}{2} \sum_{n=0}^\infty \frac{e^{-ns\left( 1+r+\frac{{\alpha'}^2}{2}\right)} }{n} e^{e^\frac{-\alpha '}{\sigma}\cdot ns} \\ &=& \frac{1}{2} \sum_{n=0}^\infty \frac{e^{-ns\left( 1+r+\frac{{\alpha'}^2}{2}-e^\frac{-\alpha '}{\sigma}\right)} }{n} = -\frac{1}{2}\ln\left(1-e^{-s\left( 1+r+\frac{{\alpha'}^2}{2}-e^\frac{-\alpha '}{\sigma}\right) }\right) \end{eqnarray*} where the last line uses that $\frac{\alpha'}{\sigma}=\frac{\alpha}{\sigma^2}\geq \alpha\cdot\left(r+\EE^0\left[e^{Z_1}\right]\right)\geq 0$ and we need to impose the condition that $e^{-\alpha\cdot\left(r+\EE^0\left[e^{Z_1}\right]\right)}\leq 1+r+\frac{\alpha^2}{2\sigma^2}$ (which, given $r>0$ and $\alpha$, will be satisfied if $\sigma>0$ is sufficiently small) to employ the identity \begin{equation}\label{geoharmo}\forall q<1\quad \sum_{n=0}^\infty\frac{q^n}{n}=\ln\frac{1}{1-q}.\end{equation} The lower bound follows simply from $$\forall n\in\NN_0\forall s>0\quad \sum_{k=0}^\infty\frac{1}{k!}\left(ns\cdot e^{ -2 \frac{\alpha'}{\sigma} - \frac{k}{\sigma^2ns}}\right)^k\geq 1$$ (for $n=0$ recall that $0^0=1$ in this paragraph by our earlier convention) as this entails (when exploiting the estimate (\ref{lowerjump1}) and finally (\ref{geoharmo}) ): \begin{eqnarray*} && \sum_{n=0}^\infty \sum_{k=0}^\infty e^{-ns}\frac{(ns)^k}{k!}\cdot \frac{e^{-rns}}{n}\PP^0\left\{B_{1}<-\frac{\alpha}{\sigma} (ns)^{\frac{1}{2}}+\frac{k}{\sigma}(ns)^{-\frac{1}{2}}\right\} \\ &\geq& \frac{1}{2\sqrt{2}} \sum_{n=0}^\infty \frac{e^{-ns\left( 1+r+{\alpha'}^2\right)}}{n} \sum_{k=0}^\infty\frac{1}{k!}\left(ns\cdot e^{-2 \frac{\alpha'}{\sigma} - \frac{k}{\sigma^2ns}}\right)^k\\ &\geq& \frac{1}{2\sqrt{2}} \sum_{n=0}^\infty \frac{e^{-ns\left( 1+r+{\alpha'}^2\right)}}{n} = -\frac{1}{2\sqrt{2}}\ln\left(1-e^{-s\left( 1+r+{\alpha'}^2\right)}\right)\end{eqnarray*} As a consequence of these estimates and using the Taylor expansion of $\exp$ around $0$, we now get the existence of two constants $c_3>0$ and $C_3>0$ (which can be computed explicitly) such that for all sufficiently small $s$, $$c_3\cdot s^\frac{1}{2}\leq \exp\left(-\sum_{n=0}^\infty\frac{e^{-rns}}{n}\PP^0\left\{X_{ns} < 0\right\} \right)\leq C_3\cdot s^\frac{1}{2\sqrt 2}.$$ Finally, we may apply identity (\ref{phase1d identity}) from Theorem \ref{phase1d} -- as this is an immediate consequence of Feller's identity \cite[p. 606, Lemma 3]{F} -- and conclude that if $G=(-\infty,\gamma)$, then $$c_3\cdot s^\frac{1}{2}\leq \rho(s)(\gamma)\leq C_3\cdot s^\frac{1}{2\sqrt 2} $$ for all sufficiently small $s>0$.

\end{ex}

\subsection{A Wiener-Hopf type result in higher dimensions}

The proof of Theorem \ref{phase1d} relies heavily on Feller's result \cite[p. 606, Lemma 3]{F} which in turn is proven by means of elementary Fourier analysis and a so-called ``basic identity'' \cite[p. 600, equation (1.9)]{F}.

Hence, if one aims at generalising Theorem \ref{phase1d} to higher dimensions, one should first of all find a multi-dimensional analogue of the identity \eqref{Feller}.

Indeed, we shall see that this is feasible. Let us for the following fix a stochastic process $X=(X_t)_{t\geq 0}$ on $\RR^d$ with stationary and independent increments.

\begin{lem} Suppose $H$ is a measurable subset of $\RR^d$, and $s>0$. Define for all $n\in\NN$ $$\forall K\in\cB\left(\RR^d\right) \quad R_n(K):=\PP^0\left[\bigcap_{1\leq i< n}\left\{X_{is}\in \complement H\right\}\cap \left\{X_{ns}\in K\cap H\right\}\right], $$ as well as $$\forall K\in\cB\left(\RR^d\right)\quad Q_n(K):=\PP^0\left[\bigcap_{1\leq i< n}\left\{X_{is}\in \complement H\right\}\cap \left\{X_{ns}\in \complement H\cap K\right\}\right] $$ (in particular $R_0=\delta_{0}\left[\cdot\cap H\right]=0$ and $Q_0=\delta_{0}\left[\cdot\cap \complement H\right]=\delta_0$). Then for all $n\in\NN_0$, $$Q_{n+1}+R_{n+1}=Q_n\ast {\PP^0}_{X_s}.$$
\end{lem}
\begin{proof} Consider a measurable $K\subseteq \RR^d$. Clearly, \begin{equation} \label{Fellergeneral} \left(Q_{n+1}+R_{n+1}\right)(K)=\PP^0\left[\bigcap_{i=1}^{n}\left\{X_{is}\in \complement H\right\}\cap \left\{X_{(n+1)s}\in K\right\}\right]. \end{equation} On the other hand, since $X$ is a Markov process, we have $$Q_n(K) =\left(P_s\left(\chi_{\complement H}\cdot\right)\right)^{\circ n}\chi_K(0)$$ (where $(P_t)_{t\geq 0}:=\left(\PP_{X_t}^0\ast\cdot\right)_{t\geq 0}$ is the translation-invariant Markov semigroup of transition functions for the process $X$ whose increments are stationary and independent), thus $$\int_{\RR^d}f(y)Q_n(dy)= \left(P_s\left(\chi_{\complement H}\cdot\right)\right)^{\circ n}f(0)$$ for all nonnegative measurable functions $f$. But this implies \begin{eqnarray*}\left(Q_{n}\ast \PP^0_{X_s}\right)(K)&=&\int_{\RR^d}\int_{\RR^d}\chi_K\left(z+y\right)\PP_{X_s}^0(dz)Q_n(dy)\\ &=& \left(P_s\left(\chi_{\complement H}\cdot\right)\right)^{\circ n}\left(\int_{\RR^d}\chi_{K-\cdot}(z)\PP_{X_s}^0(dz)\right)(0) \\ &=& \left(P_s\left(\chi_{\complement H}\cdot\right)\right)^{\circ n}\circ\left(\PP_{X_s}^0\ast\chi_K\right)(0)\\ &=& \left(P_s\left(\chi_{\complement H}\cdot\right)\right)^{\circ n}\circ P_s\chi_K(0),
\end{eqnarray*}
and the right hand side of this equation coincides with the one of identity (\ref{Fellergeneral}).
\end{proof}

Applying Fourier transforms we obtain

\begin{cor}\label{1.9general} Let us adopt the notation of the preceding Lemma and define the {\em Fourier transform} of a countable sequence $\left(\mu_n\right)_n$ of finite measures on $\RR^d$, denoted by $\widehat{\left(\mu_n\right)_n}=\widehat{\mu}:(0,1)\times\RR^d\rightarrow \CC$, by $$\forall q\in(0,1)\quad \forall\zeta \in \RR^d\quad \widehat{\left(\mu_n\right)_n}(q,\zeta)=\sum_{n=0}^\infty q^n\int_{\RR^d} e^{i\cdot{{^t}\zeta}y}\mu_n(dy)=\sum_{n=0}^\infty q^n\widehat{\mu_n}(\zeta).$$ Then for all $q\in (0,1)$, and $\zeta\in\RR^d$ the equation $$1-\widehat{\left(R_n\right)_n}(q,\zeta) =\widehat{\left(Q_n\right)_n}(q,\zeta)\left(1-q\widehat{{\PP^0}_{X_s}}(\zeta)\right)$$ holds.
\end{cor}
\begin{proof} The result of the previous Lemma reads $$\forall n\in\NN_0 \quad \widehat{Q_{n+1}}+\widehat{R_{n+1}}=\widehat{Q_n}\widehat{{\PP^0}_{X_s}}$$ when we apply the Fourier transform. After multiplication with $q^{n+1}$ and summing up over $n\in\NN_0$, one arrives at $$\forall q\in(0,1)\quad\forall \zeta\in\RR^d\quad \widehat Q(q,\zeta)- \underbrace{\chi_{\complement H}(0)}_{=\widehat{Q_0}(\zeta)}+\widehat R(q,\zeta) - \underbrace{\chi_{H}(0)}_{=\widehat{R_0}(\zeta)}= q\widehat Q(q,\zeta)\widehat{{\PP^0}_{X_s}}(\zeta) ,$$ hence $$\forall q\in(0,1)\quad \widehat R(q,\cdot)-1 = q\widehat Q(q,\cdot)\widehat{{\PP^0}_{X_s}}(\cdot) -\widehat Q(q,\cdot). $$ This is our claim.
\end{proof}

\begin{Def} A subset $A\subseteq \RR^d$ is called {\em $+$-closed} if and only if $A$ is measurable and $A+A\subseteq A$, that is sums of elements of $A$ are again elements of $A$.
\end{Def}

\begin{lem} Let $H$ be closed and convex with $0\in \partial H$. Whenever both $H$ and its complement $\complement H$ are $+$-closed, there exists a $y_H\in \RR^d$ such that $$ H=\left\{x\in\RR^d\ : \ {^t}x y_H\geq 0 \right\}.$$
\end{lem}
\begin{proof} According to the projection theorem, there exists a $y_H\in\RR^d$ such that $$ H\subseteq \left\{x\in\RR^d\ : \ {^t}x y_H\geq 0\right\}=:K.$$ We must show $H\supseteq K$. 

First, note that there is no $y\in \complement H$ such that ${^t}yy_H>0$. For, if there was one, then ${^t}(-y)y_H<0$, hence $-y\in \complement K \subseteq \complement H$, which yields, because $\complement H$ is $+$-closed, also $$0=(-y)+y\in\complement H, $$ contradicting our assumption $0\in\partial H\subseteq H$. Therefore we already have $$H\supseteq \left\{x\in\RR^d\ : \ {^t}x y_H> 0\right\}.$$ 

Now $\partial \left\{x\in\RR^d\ : \ {^t}x y_H> 0\right\} \supseteq \left\{x\in\RR^d\ : \ {^t}x y_H = 0\right\}$ (even equality holds, but this is not needed here), since whenever ${^t}x y_H = 0$, then $\frac{1}{n}y_H+x\longrightarrow x$ as $n\rightarrow \infty$ and also ${^t}\left(\frac{1}{n}y_H+x\right) y_H=\frac{1}{n}{^t}y_Hy_H>0$ for all $n\in\NN$, making $x$ a limit point of $\left\{x\in\RR^d\ : \ {^t}x y_H> 0\right\}$. 

But on the other hand, $H$ is closed by assumption and we have already proven $H\supseteq \left\{x\in\RR^d\ : \ {^t}x y_H > 0\right\}$. So we get in addition $$H\supseteq \left\{x\in\RR^d\ : \ {^t}x y_H= 0\right\}$$ and thus $H\supseteq K$ as claimed.

\end{proof}

\begin{lem}[\`a la Feller, Wiener, Hopf] \label{alaWienerHopf}Suppose $H$ is a $+$-closed set and its complement $\complement H$ is a $+$-closed set as well. Assume furthermore $0\not\in H$ (ensuring $R_0 = 0$), and let $\ln$ be the main branch of the logarithm on $\CC$. Then $$-\ln\left({1-\widehat{R}(q,\zeta)}\right)=\sum_{n=1}^\infty\frac{ q^n}{n}\int_He^{i\cdot{{^t}\zeta}x}\left({\PP^0}_{X_s}\right)^{\ast n}(dx)$$ for all $(q,\zeta)\in (0,1)\times\RR^d$ such that the left-hand side is well-defined. In general, for all $q\in (0,1)$, one has at least $${1-\widehat{R}(q,0)}=\exp\left(-\sum_{n=1}^\infty\frac{q^n\cdot\left( \widehat{\PP^0_{X_s}\left[\cdot\cap H\right]}(0)\right)^n }{n}\right).$$
\end{lem}
\begin{proof} Let $q\in (0,1)$. According to the previous Corollary \ref{1.9general}, we have \begin{equation}\forall \zeta\in U\quad \ln\frac{1}{1-q\widehat{{\PP^0}_{X_s}}(q,\zeta)}=\ln\frac{1}{1-\widehat{R}(q,\zeta)} - \ln\widehat Q(q,\zeta)\label{alaWienerHopfLogEq} \end{equation} wherever this is defined. Due to the identities $\sum_{n=1}^\infty \frac{r^n}{n}=\ln\frac{1}{1-r}$ for all $r\in B_1(0)\subset \CC$ (cf. equation (\ref{geoharmo1}) in the proof of Theorem \ref{phase1d} above) and $\widehat{{\PP^0}_{X_s}}^n=\widehat{{{\PP^0}_{X_s}}^{\ast n}}$ this can also be written as \begin{eqnarray} &&\label{alafeller}\sum_{n=1}^\infty\frac{ q^n }{n}\int_{\RR^d} e^{i\cdot{{^t}\zeta}x}\left({\PP^0}_{X_s}\right)^{\ast n}(dx) \\ \nonumber &=& \sum_{n=1}^\infty \frac{1}{n} \left(\widehat{R}(q,\zeta)\right)^n + \sum_{n=1}^\infty \frac{(-1)^n}{n}\left(\widehat{Q}(q,\zeta)-1\right)^n.\end{eqnarray} However, at least for $\zeta=0$ and arbitrary choice of $q$, one may still state identity (\ref{alaWienerHopfLogEq}) as this follows from Corollary \ref{1.9general} more or less directly: First we note that \begin{eqnarray*}&&-\ln\left(\left(1-q\widehat{{\PP^0}_{X_s}}(q,0)\right)\cdot \widehat{\left(Q_n\right)_n}(q,0)\right)\\&=&-\ln \widehat{\left(Q_n\right)_n}(q,0)-\ln\left(1-q\widehat{{\PP^0}_{X_s}}(q,0)\right)\end{eqnarray*} (as in these statements the arguments of $\ln$ are positive, hence surely in the domain of $\ln$) and written in series notation \begin{eqnarray*}&& \sum_{n=1}^\infty \frac{1}{n}\left(1-\widehat{\left(Q_n\right)_n}(q,0)\right)^n+ \sum_{n=1}^\infty\frac{q^n}{n}\widehat{{\PP^0}_{X_s}}(q,0)^n \\ &=& \sum_{n=1}^\infty \frac{(-1)^n}{n}\left(\left(1-q\widehat{{\PP^0}_{X_s}}(q,0)\right)\cdot \widehat{\left(Q_n\right)_n}(q,0)-1\right)^n.\end{eqnarray*} But Corollary \ref{1.9general} implies $$ \forall n\in\NN\quad\frac{(-1)^n}{n}\left(\left(1-q\widehat{{\PP^0}_{X_s}}(q,0)\right)\cdot \widehat{\left(Q_n\right)_n}(q,0)-1\right)^n = \frac{1}{n}\widehat{\left(R_n\right)_n}(q,0).$$ Combining these two equations yields (\ref{alaWienerHopfLogEq}). Next, note that $$\mu_{R,q}:=\sum_{n=0}^\infty q^nR_n$$ is still a finite measure -- concentrated on $H$ -- and thus possesses a Fourier transform. Analogously, the measure $\mu_{Q,q}:=\sum_{n=0}^\infty q^nQ_n$ is concentrated on $\complement H$ and also has a Fourier transform as it is finite. Now, for arbitrary $n\in\NN$, the properties of the Fourier transform imply \begin{eqnarray*}\left(\widehat{R}(q,\cdot)\right)^n&=& \left(\widehat{\mu_{R,q}}\right)^n=\widehat{{\mu_{R,q}}^{\ast n}}, \\ \left(\widehat{Q}(q,\cdot)-1\right)^n &=& \left(\widehat{\mu_{Q,q}-\delta_0}\right)^n=\widehat{\left( \left({\mu_{Q,q}-\delta_0}\right)^{\ast n} \right)}.\end{eqnarray*} But since $H$ and $\complement H$ are $+$-closed sets, ie $H+H\subseteq H$ and $\complement H+\complement H\subseteq\complement H$ , the measures on the right hand sides of these two equations, ${{\mu_{R,q}}^{\ast n}}$ and $\left({\mu_{Q,q}-\delta_0}\right)^{\ast n}$, have to be (signed) measures on $H$ and $\complement H$, respectively. Let us now split the sum in (\ref{alafeller}) and insert the terms we have previously identified:
\begin{eqnarray} \nonumber&& \forall \zeta\in U \\ &&\sum_{n=1}^\infty\frac{ q^n }{n}\int_{H} e^{i\cdot{{^t}\zeta}x}\left({\PP^0}_{X_s}\right)^{\ast n}(dx) + \sum_{n=1}^\infty\frac{ q^n }{n}\int_{\complement H} e^{i\cdot{{^t}\zeta}x}\left({\PP^0}_{X_s}\right)^{\ast n}(dx) \nonumber \\ &=& \sum_{n=1}^\infty \frac{1}{n} \widehat{\left({\mu_{R,q}}^{\ast n}\right)}(\zeta) + \sum_{n=1}^\infty \frac{(-1)^n}{n}\widehat{\left({\mu_{Q,q}-\delta_0}\right)^{\ast n}}(\zeta).\end{eqnarray} It is the injectivity of the Fourier transform that yields from this \begin{eqnarray} \nonumber &&\sum_{n=1}^\infty\frac{ q^n }{n}\left(\left({\PP^0}_{X_s}\right)^{\ast n}\left(\cdot\cap H\right)\right) + \sum_{n=1}^\infty\frac{ q^n }{n} \left(\left({\PP^0}_{X_s}\right)^{\ast n}\left(\cdot\cap \complement H\right)\right) \\ \nonumber &=& \sum_{n=1}^\infty \frac{1}{n}\left({\mu_{R,q}}^{\ast n}\right) + \sum_{n=1}^\infty \frac{(-1)^n}{n}\left({\mu_{Q,q}-\delta_0}\right)^{\ast n}.\end{eqnarray} Either side of this equation equals the sum of two (signed measures), and we recall that the first measure on the left hand side and first measure on the right hand side are both concentrated on $H$, whilst the second measure on the left hand side as well as the second measure on the right hand side are both concentrated on $\complement H$. The only way for this to be true is that the two measures that are concentrated on each of $H$ or $\complement H$ are equal: $$\sum_{n=1}^\infty\frac{ q^n }{n}\left(\left({\PP^0}_{X_s}\right)^{\ast n}\left(\cdot\cap H\right)\right) = \sum_{n=1}^\infty \frac{1}{n}\left({\mu_{R,q}}^{\ast n}\right),$$ and also $$ \sum_{n=1}^\infty\frac{ q^n }{n} \left(\left({\PP^0}_{X_s}\right)^{\ast n}\left(\cdot\cap \complement H\right)\right) = \sum_{n=1}^\infty \frac{(-1)^n}{n}\left({\mu_{Q,q}-\delta_0}\right)^{\ast n} ,$$ the former identity being exactly what the statement of the Lemma expresses in the language of Fourier transforms.
\end{proof}

\subsection{Continuity corrections in higher dimensions}

On the basis of Lemma \ref{alaWienerHopf}, we may now partially generalise Theorem \ref{phase1d} to higher dimensions when we require $G$ (the set that we refer to the exercise region) to be $+$-closed set.

\begin{Th}\label{phase} Let $X$ be a $d$-dimensional L\'evy basket with an associated family of risk-neutral measures $\PP^\cdot$ and discount rate $r>0$. Also, consider a $G$ of the shape $G=\gamma+H$ for some $\gamma\in\RR^d$ and some measurable set $H$ such that $0\not \in H$, and such that both $H$ as well as $\complement H$ are $+$-closed.
Then for all $s>0$, \begin{eqnarray*} \rho(s)(\gamma)&=& \left. \sum_{n=1}^\infty e^{-rnt} \PP^{\gamma} \left[ \bigcap_{i=1}^{n-1}\left\{ X_{it}\in\complement G\right\}\cap\left\{X_{nt}\in G\right\}\right]\right|^{t\downarrow 0}_{t=s} \\&=& \left. \sum_{n=1}^\infty e^{-rnt} \PP^{0} \left[ \bigcap_{i=1}^{n-1}\left\{ X_{it}\in\complement H\right\}\cap\left\{X_{nt}\in H\right\}\right]\right|^{t\downarrow 0}_{t=s}\\ &=& \left. \exp\left(-\sum_{n=1}^\infty \frac{e^{-rns}}{n}\PP^{0} \left\{X_{ns}\in G-\gamma \right\} \right)\right|^{t=s}_{t\downarrow 0} \end{eqnarray*} (and the limit on the right exists whenever $\rho(s)(\gamma)$ is well-defined).
\end{Th}
\begin{proof} The Theorem follows directly from Lemma \ref{alaWienerHopf} --- in the same manner in which Theorem \ref{phase1d} followed from Feller's original result \cite[p. 606, Lemma 3]{F}: For, the second equation in Lemma \ref{alaWienerHopf} may be read $${1-\sum_{n=0^\infty}q^nR_n\left[\RR^d\right] } = \exp\left(-\sum_{n=1}^\infty\frac{q^n\cdot\left( \PP^0_{X_s}\left[H\right]\right)^n }{n}\right)$$ that is \begin{eqnarray*} && 1-\sum_{n=0^\infty}q^n\PP^0\left[\bigcap_{1\leq i< n}\left\{X_{is}\in \complement H\right\}\cap \left\{X_{ns}\in H\right\}\right]\\ & = & \exp\left(-\sum_{n=1}^\infty\frac{q^n\cdot\left( \PP^0 \left\{X_s\in H\right\}\right)^n }{n}\right) \end{eqnarray*} for all $q\in (0,1)$, in particular for $q=e^{-rs}$.

\end{proof}

\begin{cor}\label{phasecorollary} Let $X$ be a $d$-dimensional Black-Scholes model with discount rate $r>0$ and volatility vector $\sigma$ (i.e. $\left(X_i\right)_t= \left(X_i\right)_0+ \sigma_iB^{(i)}_t+\left(r-\frac{\sigma_i^2}{2}\right)t$ for all $t>0$ and independent normalised Brownian motions $B^{(i)}$, $i\in\{1,\dots,d\}$). Let $\PP^\cdot$ denote the corresponding family of risk-neutral measures. Also, consider a region $G$ of the shape $G=\gamma+\left\{x\in\RR^d \ : \ {^t}\alpha x< 0\right\}$ for some $\alpha,\gamma\in\RR^d$. Let $\delta:= \sqrt{\sum_{i=1}^d\left|\alpha_i\sigma_i\right|^2}$ and $\mu:={^t}\left(r-\frac{{\sigma_1}^2}{2}, \dots, r-\frac{{\sigma_d}^2}{2}\right)$. 

Then, whenever ${^t}\alpha\mu\geq 0$, there exist constants $c_0, C_0>0$ such that for all sufficiently small $s>0$, \begin{eqnarray*} c_0s^\frac{1}{2}\leq \rho(s)(\gamma) \leq C_0s^\frac{1}{2\sqrt{2}}\end{eqnarray*} If both ${^t}\alpha\mu\leq 0$ and $r> \frac{\left|{^t}\alpha\mu\right|^2}{2\delta^2}$, there exist constants $c_1, C_1>0$ such that for all sufficiently small $s>0$, \begin{eqnarray*}c_1s^\frac{1}{\sqrt{2}} \leq \rho(s)(\gamma)\leq C_1s^\frac{1}{{2}}.\end{eqnarray*} 

\end{cor}

\begin{proof}[Proof of Corollary \ref{phasecorollary}] Observe that \begin{eqnarray*}\PP^{\gamma} \left\{X_{ns}\not\in G \right\} &=&\PP^{0} \left\{X_{ns}\not\in (G-\gamma) \right\} =\PP^{0} \left\{{^ t}\alpha X_{ns}\geq 0 \right\}\end{eqnarray*}
But since the components of $X$ are independent multiples of Brownian motions with linear drift at rates $\mu_1,\dots,\mu_d$, respectively, the process $\left({^ t}\alpha X_t\right)_{t\geq 0}$ is a multiple of a normalised Brownian motion with linear drift at rate ${^t}\alpha\mu$: For each $t>0$, the random variable $${^ t}\alpha X_t=\sum_{i=1}^d\left(\alpha_i\sigma_iB^{(i)}_t\right)-\sum_{i=1}^d\alpha_i\left(r-\frac{{\sigma_i}^2}{2}\right)t$$ is distributed according to \begin{equation}\PP^0_{{^t}\alpha X_t}=\nu_{{^t}\alpha\mu\cdot t, \sum_{i=1}^d\left|\alpha_i\sigma_i\right|^2t} = \nu_{{^t}\alpha\mu\cdot t, \delta^2t}\label{alpha X distribution} .\end{equation} 

Hence for any normalised $1$-dimensional Brownian motion $B$, one has \begin{eqnarray*}\PP^{\gamma} \left\{X_{ns}\not\in G \right\} &=&\PP^{0} \left\{\delta\cdot B_{ns}+{^t}\alpha\mu\geq 0\right\}\end{eqnarray*} and \begin{eqnarray*}\PP^{\gamma} \left\{X_{ns}\in G \right\} &=&\PP^{0} \left\{\delta\cdot B_{ns}+{^t}\alpha\mu< 0\right\}. \end{eqnarray*}

We can now apply the estimates from the proof of Corollary \ref{BlackScholes1d} to the process $\left(\delta B_t +{^t}\alpha\mu t\right)_{t\geq 0}$ {\em en lieu} of what is $X$ there. 

For instance, if ${^t}\alpha\mu\geq 0$, we shall employ the bounds $$\frac{e^{-\frac{\left|{^t}\alpha\mu\right|^2}{\delta^2} ns}}{2\sqrt{2}}\leq \PP^\gamma\left\{ X_{ns}\in G \right\}\leq \frac{e^{-\frac{\left|{^t}\alpha\mu\right|^2}{2\delta^2}ns}}{2}$$ and if ${^t}\alpha\mu\leq 0$, we will use the estimates 
$$ \frac{e^{-\frac{\left|{^t}\alpha\mu\right|^2}{2\delta^2}ns} }{2} \leq \PP^\gamma\left\{ X_{ns}\in G \right\}\leq \frac{ e^{\frac{\left|{^t}\alpha\mu\right|^2}{2\delta^2}ns} }{\sqrt{2}}.$$
\end{proof}

\pagebreak 
\appendix{\Large\bf Appendix \bigskip \\ Reformulation of the perpetual Bermudan pricing problem in $L^1$ and $L^2$}

\section{Non-applicability of the \\ $L^2(\RR^d)$ Spectral Theorem}

Consider a $d$-dimensional L\'evy basket $X$ with associated family of risk-neutral probability measures $\PP^\cdot$ and discount rate $r>0$.

Fixing $h>0$ and defining $$P:=\pi_{L^2\left(\RR^d\setminus G\right)},\quad A:=A^h:= \II - 
e^{-rh}\PP_{X_0-X_h}\ast\cdot=\left(\delta_0 - e^{-rh}\PP_{X_0-X_h}\right)\ast\cdot,$$ we can rewrite the result of Lemma \ref{convolveproject} as follows:
\begin{equation}\label{1.1} PA\left(V_G^h-g_1\right)=-PAg_1\end{equation} where we assume that $g$ has a square-integrable extension from $G$ to the whole of $\RR^d$; given this assumption, the $g_1\in L^2(\RR^d)$ of the previous identity can be any such extension.

We will suppress the superscript of $A$ for the rest of this paragraph.

Also, without loss of generality, we will assume in this Chapter that the components of the basket $X$ when following the Black-Scholes model all have volatility $1$.

\begin{lem} Let $X$ be a L\'evy basket with associated family of risk-neutral probability measures $\PP^\cdot$ and discount rate $r>0$. Then $A$ and $PA\restriction L^2(\RR^d\setminus G)$ are invertible. Furthermore, the $L^2$ norm of $A$ is bounded by $\left(1+e^{-rh}\right)^\frac{1}{2}$, if $X=\mu\cdot+B$ (thus $\mu=\left(r-\frac{1}{2}\right)_{i=1}^d$) where $B$ is a standard Brownian motion. Moreover, $A$ is a contraction if $\mu=0$. 
\end{lem}
\begin{proof} Suppose $0\neq u\in L^2(\RR^d\setminus G)$ and $u$ is bounded. Then $\esssup |u|\neq 0$ and we may choose a set $H\subset G$ of positive Lebesgue measure such that $e^{-rh}\esssup|u| <|u(x)|$ for all $x\in H$ (this is possible because $r,h>0$ and therefore $e^{-rh}<1$), we deduce \begin{eqnarray*}\forall x\in H\quad \left|e^{-rh}\PP_{X_0-X_h}\ast u(x)\right| &\leq& e^{-rh} \esssup|u|\\&<& u(x), \end{eqnarray*} which means that $$\forall x\in H \quad PAu(x)=u(x)-e^{-rh}\PP_{X_0-X_h}\ast u(x)\neq 0,$$ hence $PAu\neq 0$ (for $H$ has positive Lebesgue measure). So $$\ker PA\restriction L^2(\RR^d\setminus G)=\ker PA\cap L^2(\RR^d\setminus G)=\{0\}$$ and we are done for the invertibility of $PA\restriction L^2(\RR^d\setminus G)$. Similarly, one can prove the invertibility of $A$. Finally, $A$ is seen to be a contraction by application of the Fourier transform: The Fourier transform is an $L^2$ isometry (by Plancherel's Theorem), thus \begin{eqnarray*}&&\left\|\left(\delta_0 - e^{-rh}\PP_{X_0-X_h}\right)\ast f\right\|_{L^2\left(\RR^d\right)}\\ &=&\left\|\left(\left(\delta_0 - e^{-rh}\PP_{X_0-X_h}\right)\ast f\right)^{\widehat{}}\right\|_{L^2\left(\RR^d\right)} \\ &=& \left\|\left(1 - e^{-rh}\widehat{\PP_{X_0-X_h}}\right)\cdot \widehat f\right\|_{L^2\left(\RR^d\right)} \\ &=& \left\|\left(1 - e^{-rh}e^{ih{^t}\mu\cdot}e^{-|\cdot|^2h/2}\right) \widehat f\right\|_{L^2\left(\RR^d\right)}\\ &=& \left(\int_{\RR^d}\left|1-e^{-rh+ih{^t}\mu x- \frac{|x|^2h}{2}}\right|^2\left|\widehat f(x)\right|^2 \ dx\right)^\frac{1}{2}. \end{eqnarray*}Now, the factor in front of $\left|\widehat f(x)\right|^2$ in the last line can be bounded by $\left(1+e^{-rh}\right)^{2}$, and it is strictly less than one for $\mu=0$. Using Plancherel's Theorem again, this yields the result. 
\end{proof}

Now, this is sufficient to apply a Wiener-Hopf factorisation (for a general treatment of this kind of factorisations, one may consult e.g. Speck \cite{Sp85}, our application uses in particular \cite[1.1, Theorem 1]{Sp85}) and state

\begin{Th} \label{l2whf} Let $G\subseteq\RR^d$ and let $X$ be a L\'evy basket with associated family of risk-neutral probability measures $\PP^\cdot$ and discount rate $r>0$. Then $V_G^h$, the expected payoff of a perpetual Bermudan option for $G$ with exercise mesh size $h>0$ and payoff function $g$, is -- using the above notation -- given by $$V_G^h=g_1-\left(PA\restriction L^2(\RR^d\setminus G)\right)^{-1}PAg_1=g_1-A_+^{-1}PA_-^{-1}PAg_1$$ where $A=A_-A_+$ is a Wiener-Hopf factorisation of $A$.
\end{Th}

We observe the following:

\begin{lem} The Hilbert space operator $A:L^2(\RR^d,\CC)\rightarrow L^2(\RR^d,\CC)$ is normal.
\end{lem}
\begin{proof} We define $p:=e^{-rh} \frac{d\PP_{X_0-X_h}}{d\lambda^d}$ (where $\lambda^d$ is the $d$-dimensional Lebesgue measure) and via the Fubini Theorem one has for every $f,g\in L^2(\RR^d)$
\begin{eqnarray*} \langle Af,g\rangle&=& \int_{\RR^d}\int_{\RR^d}\left(\delta_0-p\right) (x-y)f(y)dy \bar g(x) dx \\ &=& \int_{\RR^d}\int_{\RR^d}f(y)\left(\delta_0-p\right) (x-y) \bar g(x) dx dy\\ &=& \langle f, \overline{\left(\delta - \bar p\circ(-\II)\ast g\right) }\rangle,\end{eqnarray*} that is $$A^*=\left(\delta_0-\bar p\circ(-\II)\right)\ast\cdot $$ But since the convolution is associative and commutative, this implies \begin{eqnarray*}A^*A&=& \left(\delta_0- \bar p\circ(-\II)\right)\ast\left(\delta_0-p\right)\ast\cdot \\ &=& \left(\delta_0-p\right)\ast \left(\delta_0-\bar p\circ(-\II)\right)\ast\cdot \\&=& AA^*.\end{eqnarray*}
\end{proof}

However, it will not be possible to find a basic system of eigenvectors and eigenvalues for this operator, since 

\begin{lem} The operator $A$ fails to be compact.
\end{lem} 
\begin{proof} Any normalised basis provides a counterexample for the compactness assertion.
\end{proof}

Therefore, the equation (\ref{1.1}) cannot easily be applied to compute the expected option payoff by means of a spectral analysis. Thus, our examination of the Hilbert space approach in the second part of this Chapter has led to a negative outcome.

However, one can also conceive of the operators $A^h$ as operators on the Banach space $L^1\left(\RR^d\right)$:

\section{The $L^1$ operator equation: \\ analyticity in the exercise mesh size}

From now on, $h$ will no longer be fixed and we will therefore write $A^h$ instead of $A$.

If we now assume $g_1$ to be an integrable extension of $G$ to the complement of $\complement G$ as an element of Quite similarly to \ref{l2whf}, we can prove

\begin{Th} Let $G\subseteq\RR^d$ and let $X$ be a L\'evy basket with associated family of risk-neutral probability measures $\PP^\cdot$ and discount rate $r>0$. Then $V_G^h$, the expected payoff of a perpetual Bermudan option for $G$ with exercise mesh size $h>0$ and payoff function $g$, is -- using the above notation -- given by $$V_G^h=g_1-\left(PA\restriction L^1(\RR^d\setminus G)\right)^{-1}PAg_1=g_1-A_+^{-1}PA_-^{-1}PAg_1,$$ where $A=A_-A_+$ is a Wiener-Hopf factorisation of $A$.
\end{Th}

It suffices to observe that $A^h$ is -- due to the $L^1$ norm estimate for the convolution of two integrable functions (as the product of the norms of the convolved functions) -- also a bounded operator on $L^1\left(\RR^d\right)$.

We shall now identify $g$ and $g_1$.

\begin{Th} With the notation previously introduced, we define $E$ to be the semigroup $$(E_t)_{t\geq 0}=\left(e^{-rt}\nu_{\mu t,t}\ast\cdot\right)_{t\geq 0}=\left( e^{-rt}g_{\mu t,t}\ast\cdot\right)_{t\geq 0}\in L\left(L^1(\RR^d),L^1(\RR^d)\right)^{[0,+\infty)},$$ where $$g_{\mu t,t}=(2\pi t)^{-d/2}e^{-|\mu t -\cdot|^2/(2t)}$$ is the distribution of the logarithmic price vector at time $t$. Suppose $x\not\in G$ and, with the notation from the previous chapters, $X$ is a (normalised) Brownian motion with (possibly zero) drift ({\em Black-Scholes model}). Then $h\mapsto V_G^h(x)$ is real analytic in $h$ on $(0,+\infty)$ as function with range in the Banach space $L^1(\RR^d)$. 
\end{Th}

\begin{proof} It is obvious that $E$ is a semigroup. According to \cite[Theorem 1.48]{D}, the set $$\cE:=\left\{f\in L^1(\RR^d) \ : \ t\mapsto E_tf \in L^1(\RR^d)\text{ entire} \right\}$$ is dense in $L^1(\RR^d)$. Hence it is possible to approximate every $g$ by a sequence $\{g_k\}_k\subset\cE$ in $L^1(\RR^d)$. Since $$\forall t\geq 0 \quad \left\|E_t\right\|_{L^1(\RR^d)}\leq e^{-rt}\leq 1,$$ we obtain $E_tg_k\rightarrow E_tg$ for $k\rightarrow\infty$ uniformly in $t$ on $\RR_+$, where $E_tg_k$ is entire for every $t>0$ and $k\in\NN$. Thus, $t\mapsto E_tg$, and thereby $t\mapsto Pg-PE_tg$, is an analytic function on $(0,+\infty)$ taking values in the Banach space $L^1(\RR^d)$. 
Now observe that for arbitrary open $U\subset\subset\RR_+$ (the symbol ``$\subset\subset$'' indicating that $U$ is contained in a compact subset of $\RR_{>0}$) the following equations hold: \begin{eqnarray}\nonumber \forall t\in U \quad V_G^t&=&g-\left(\II-PE_t\restriction PL^1(\RR^d)\right)^{-1}\left(Pg-PE_tg\right) \\ \nonumber &=& g-\left(\sum_{k=0}^\infty(PE_t)^k\right)\left(Pg-PE_tg\right) \\ \nonumber &=& g -\sum_{k=0}^\infty (PE_t)^{k}Pg +\sum_{k=0}^\infty (PE_t)^{k+1}g \\ \label{harmonic} &=& \sum_{k=0}^\infty (PE_t)^{k}\left(g-Pg\right),\end{eqnarray} since the sums converge uniformly in $t$ on $U\subset\subset \RR_+$, yielding the analyticity of $t\mapsto V_G^t$ as a function whose range lies in the Banach space $L^1(\RR^d)$.
\end{proof}

\begin{lem} \label{deriv} Let $u>0$, $n\in\NN$. Then the equation $$\frac{d^n}{du^n}E_u = \left(-r+\frac{1}{2}\Delta-{{^t}\mu} \nabla\right)^n g_{t\mu ,t}\ast \cdot $$ holds (where ${^t}y$ denotes the transpose of a vector $y$). In particular, if $\frac{1}{2}\Delta f-\mu\nabla f= \lambda f$ for some $\lambda, f$, $$\frac{d^n}{du^n}E_u f= \left(-r+\lambda\right)^n g_{t\mu ,t}\ast f.$$ 
\end{lem}
\begin{proof} According to Davies \cite[Proof of Theorem 2.39]{D}, we have \begin{equation} \label{infgen} \frac{d^n}{du^n}E_u=\left(ZE_{u/n}\right)^n,\end{equation} where $Z$ denotes the infinitesimal generator of the semigroup $E$. Now, define $C$ to be the convolution operator semigroup $\left(g_{t\mu,t}\ast\cdot\right)_{t\geq 0}$ of (normalised) Brownian motion with drift $\mu$ (as before denoting by $g_{z,\sigma^2}$ the Lebesgue density of the Gaussian distribution centered around $z$ of variance $\sigma^2$ for all $z\in\RR^d$ and $\sigma>0$). It is well-known (cf. e.g. \cite[p. 352]{RY}) that the infinitesimal generator of this semigroup $C$ is $$L:=\frac{1}{2}\Delta +{{^t}\mu} \nabla.$$ By our requirements on $f$, $Lf=0$ on $U$. Furthermore, $L$ and $C$ commute: $$\forall t\geq 0 \quad C_tL=LC_t.$$ 
Thus, \begin{eqnarray*} \forall t\geq 0 \quad ZE_t &=& \frac{d}{d}E_t =\frac{d}{dt}\left(e^{-rt}\cdot C_t\right)\\ &=& -re^{-rt}C_t+e^{-rt}\frac{d}{dt}C_t \\&=& e^{-rt}C_t(-r+L) ,\end{eqnarray*} which due to equation (\ref{infgen}) already suffices for the proof of the Lemma in the general case. And if $f$ is an eigenfunction of $L$ for the eigenvalue $\lambda$, one has $(-r+L)^nf=(-r+\lambda)^nf$.

\end{proof}

\begin{Th}\label{taylor} The Taylor series for the expected payoff of a perpetual Bermudan option as a function of the exercise mesh with respect to a fixed exercise region $G$ is for all $s>0$: \begin{eqnarray*} \forall t>0 \quad V_G^t &=& \sum_{k=0}^\infty (t-s)^k \sum_{m=1}^\infty e^{-rms}\sum_{\begin{array}{c} l_1+\dots+l_m=n \\ (l_1,\dots,l_m)\in{\NN_0}^m \end{array} } \left(\prod_{i=1}^m\frac{1}{l_i !}\right) \\ && \left( \chi_{\RR^d\setminus G}\cdot \left(g_{s\mu,s}\ast \cdot \right)\left(-r+\frac{1}{2}\Delta+{{^t}\mu}\nabla\right)^{\circ l_i}\right)\left(\chi_Gg\right),\end{eqnarray*} where, in order to avoid confusion with pointwise exponentiation, $A^{\circ k}$ denotes $A^k$ for any operator $A$.
\end{Th}
\begin{proof} We know about the real analyticity of $t\mapsto E_t$ on $\RR^{>0}$ and even, thanks to the previous Lemma, the explicit Taylor series. Thereby we also have the Taylor series for $t\mapsto PE_t$. So we can use equation (\ref{harmonic}) and see by means of a binomial expansion
\begin{eqnarray*}\sum_{k=0}^\infty \left(PE_t\right)^k &=&\sum_{k=0}^\infty \left(P\sum_{\ell= 0}^\infty \frac{(t-s)^\ell}{\ell !}\left(e^{-rs}\left(-r+L\right)\right)^\ell C_s\right)^k\\ &=& \sum_{n=0}^\infty \sum_{m=1}^\infty \sum_{\begin{array}{c} l_1+\dots+l_m=n \\ (l_1,\dots,l_m)\in{\NN_0}^m \end{array} } \prod_{i=1}^m \frac{(t-s)^{l_i}}{l_i!}e^{-rs}P(-r+L)^{l_i}C_s \\&=& \sum_{n=0}^\infty (t-s)^n \sum_{m=1}^\infty e^{-rms} \\&& \sum_{\begin{array}{c} l_1+\dots+l_m=n \\ (l_1,\dots,l_m)\in{\NN_0}^m \end{array} } \prod_{i=1}^m \frac{1}{l_i!} P(-r+L)^{l_i}C_s.\end{eqnarray*} 
\end{proof}

This Taylor series fails to provide any straightforward possibility for the computation of $V_G$. Instead we state the following immediate Corollary of equation (\ref{harmonic}):

\begin{cor}\label{firstderiv} With the notation as in the previous Theorem, \begin{eqnarray*} \forall s>0\quad \frac{d}{ds}V_G^s &=& \frac{d}{ds}\sum_{m=1}^\infty e^{-rms}\left( \chi_{\RR^d\setminus G}\cdot \left(g_{s\mu,s}\ast \cdot \right)\right)^{\circ m} \left(\chi_G\cdot g\right).\end{eqnarray*} 
\end{cor}

{\bf Acknowledgements.} The author would like to thank the German Academic Exchange Service for the pre-doctoral research grant he received ({\em Doktorandenstipendium des Deutschen Akademischen Austauschdienstes}) and the German National Academic Foundation ({\em Studienstiftung des deutschen Volkes}) for their generous support in both financial and non-material terms. He also owes a huge debt of gratitude to his supervisor, Professor Terry J Lyons, as well as to Dr Ben Hambly and Professor Alexander Schied for their constructive comments on a previous version of this paper.

\pagebreak

\end{document}